\documentclass[12pt]{article}
\usepackage{array}

\usepackage{bm}
\newcommand{\cc}{\bm{C}}
\newcommand{\dd}{\bm{D}}
\newcommand{\ggg}{\bm{g}}
\newcommand{\ggt}{\bm{G}'}
\newcommand{\g}{\bm{G}}
\newcommand{\hh}{\bm{H}}
\newcommand{\ii}{\bm{I}}
\newcommand{\lll}{\bm{l}}
\newcommand{\Lll}{\bm{L}}
\newcommand{\Lllh}{\bm{L}^{1/2}}
\newcommand{\ddh}{{\dd}^{1/2}}
\newcommand{\qq}{\bm{q}}
\newcommand{\uu}{\bm{u}}
\newcommand{\sss}{\bm{S}}
\newcommand{\ww}{\bm{W}}

\newcommand{\zz}{\bm{0}}
\newcommand{\qqq}{\bm{Q}}
\newcommand{\uuu}{\bm{U}}

\newcommand{\ggamma}{\bm{\Gamma}}

\newcommand{\ggammat}{\bm{\Gamma}'}
\newcommand{\Llambda}{\bm{\Lambda}}
\newcommand{\Llambdaih}{\bm{\Llambda}^{-1/2}}
\newcommand{\Llambdain}{\Llambda^{-1}}
\newcommand{\llambda}{\bm{\lambda}}
\newcommand{\s}{\bm{\Sigma}}
\newcommand{\si}{\bm{\Sigma}^{-1}}

\newcommand{\gtilde}{\widetilde{\bm{G}}}
\newcommand{\gtildet}{\widetilde{\bm{G}'}}
\newcommand{\iip}{\bm{I}_p}
\newcommand{\sh}{\widehat{\bm{\Sigma}}}

\newcommand{\xxi}{\bm{\xi}}
\newcommand{\xxxi}{\bm{\Xi}}
\newcommand{\xxxii}{\xxxi^{-1}}
\newcommand{\xxxih}{\widehat{\bm{\Xi}}}
\newcommand{\ddd}{\bm{d}}

\newcommand{\op}{{\cal O}(p)}
\newcommand{\opp}{{\cal O}^+(p)}

\newcommand{\dts}{\prod_{i=1}^p \lambda_i^{-\frac{n}{2}}}
\newcommand{\lff}{\prodp l_i^{\frac{n-p-1}{2}}}
\newcommand{\lfs}{\prod_{j<i} (l_j-l_i)}
\newcommand{\prodp}{\prod_{i=1}^p}
\newcommand{\tr}{\mathop{\rm tr}}
\newcommand{\etr}{\mathop{\rm etr}}
\newcommand{\diag}{\mathop{\rm diag}}
\newcommand{\proof}{\noindent{\bf Proof.}\quad}
\newcommand{\qed}{\hbox{\rule[-2pt]{3pt}{6pt}}}

\newcommand{\sumi}{\sum_{i=1}^p}
\newcommand{\sumj}{\sum_{j=1}^p}
\newcommand{\plim}{\stackrel{p}{\rightarrow}}
\newcommand{\dlim}{\stackrel{d}{\rightarrow}}

\newtheorem{cor}{Corollary}
  
\newtheorem{lem}{Lemma}
\newtheorem{theorem}{Theorem}

\renewcommand{\aa}{\alpha}
\newcommand{\bb}{\beta}
\newcommand{\bdiva}{\bb/\aa}
\newcommand{\tn}{u}
\newcommand{\mdiff}{{\bar m}}
\newcommand{\covernum}{\tau}
\newcommand{\covertotal}{T}
\newcommand{\lossd}{d}
\newcommand{\os}{O^{(\covernum)}}

\title{Asymptotic Distribution of Wishart Matrix for Block-wise Dispersion of Population Eigenvalues}
\author{Yo Sheena\thanks{Department of Economics, Shinshu University}
  \ and Akimichi Takemura\thanks{Graduate School of Information Science and Technology, University of Tokyo}}
\date{September, 2006}

\begin{document}
\maketitle
\begin{abstract}
  This paper deals with the asymptotic distribution of Wishart matrix
  and its application to the estimation of the population matrix
  parameter when the population eigenvalues are block-wise infinitely
  dispersed. We show that the appropriately normalized eigenvectors
  and eigenvalues asymptotically generate two Wishart matrices and one
  normally distributed random matrix, which are mutually independent.
  For a family of orthogonally equivariant estimators, we calculate
  the asymptotic risks with respect to the  entropy or the quadratic loss
  function and derive the asymptotically best estimator among the
  family. We numerically show 1) the convergence in both the
  distributions and the risks are quick enough for a practical use, 2)
  the asymptotically best estimator is robust against the deviation of 
  the population eigenvalues from the block-wise infinite dispersion.
\end{abstract}

\section{Introduction}
\label{sec:intro}
Suppose that a $p$-dimensional random vector $\bm{y}$ has the
covariance matrix $\s$. 
The inference for $\s$ has been studied in enormous
amount of literature and is still an important topic from both
theoretical and practical points of view.  Often we assume some
structure of $\s$, i.e., restriction on its parameter space 
$\{\s \mid \s>0\}$.
A structure, in some cases, arises from a theoretical
reason behind the data. In other cases, it appears as a result of
exploratory analysis such as principle component analysis or
exploratory factor analysis.

For example suppose that 
$\bm{y}$ is generated in the following multivariate linear model;
\begin{equation}
\label{basic_model}
\bm{y}=\bm{B}\bm{x}+\bm{e},
\end{equation}
where $\bm{B}$ is a $p\times m$ coefficient (factor loading) matrix
with $\mbox{rank}\;{\bm B}=m$, $\bm{x}$ is a latent $m\times 1$ random
vector (common factor) and $p\times 1$ vector $\bm{e}$ is an error
term (unique factor) which is independently distributed from $\bm{x}$.
If we further assume that $\bm{e}$ has $\sigma^2\ii_p$ ($\ii_p$:
$p$-dimensional identity matrix) as its covariance matrix, $\s$ is
written as
$$
\s=\bm{B}\s_x\bm{B}'+\sigma^2\ii_p,
$$
where $\s_x$ is the nonsingular covariance matrix of $\bm{x}$. In this
case $\s$ has the eigenvalues $\lambda_1 \ge\cdots\ge \lambda_p$ given by
\begin{equation}
\label{structure_S}
\lambda_i=\left\{
\begin{array}{ll}
\tau_i+\sigma^2,&\mbox{ if }i=1,\ldots,m,\\
\sigma^2,&\mbox{ if }i=m+1,\ldots,p.
\end{array}
\right.
\end{equation}
where $\tau_i>0$, $i=1,\ldots,m$,  are the  eigenvalues of
$\bm{B}\s_x\bm{B}'$. It is often observed that $\sigma^2$ is quite
small compared to $\tau_i$'s, which means that the first group of
eigenvalues $(\lambda_1,\ldots, \lambda_m)$ is very
large compared to the second group 
$(\lambda_{m+1},\ldots,\lambda_p)$. In this paper we call this
state as ``(two-)block-wise dispersion'' of the population
eigenvalues.

What would happen to the sample covariance matrix, 
when the eigenvalues of population covariance matrix are
``infinitely'' dispersed?
This is an interesting
question from a theoretical standpoint. Takemura and Sheena (2005) and
Sheena and Takemura (2006) deal with this problem 
under ``total dispersion'' of population
eigenvalues, namely
$$
(\lambda_2/\lambda_1,\lambda_3/\lambda_2, \ldots, \lambda_p/\lambda_{p-1})\rightarrow \bm{0}.
$$
This paper is a generalization of Takemura and Sheena (2005)
from a theoretical point of view, while the practical motivation is as
follows; as we saw above, we often come across a practical situation
where the population eigenvalues are block-wise dispersed. 
It is helpful for the inference on $\s$  in practical situations to 
understand the behavior of the sample covariance matrix, when the population
eigenvalues are block-wise ``infinitely'' dispersed.
The state of the population eigenvalues being infinitely dispersed
is a theoretical approximation, 
but understanding the 
limiting behavior leads to a better insight on its neighborhood where
the eigenvalues are ``largely'' dispersed.

Now we formally state the framework of this paper. Let $\sss=(s_{ij})$ be distributed according to Wishart distribution
$\ww_p(n,\s)$, where $p$ is the dimension, $n$ is the degrees of
freedom, and $\s$ is the covariance matrix. 
The spectral decompositions of $\s$ and $\sss$  are given by
$$
\s=\ggamma\Llambda\ggammat, \qquad   \sss=\g\Lll\ggt,
$$
where $\g, \ggamma \in \op,$ the group of $p\times p$ orthogonal matrices, and 
$\Llambda=\diag(\lambda_1,\ldots,\lambda_p)$,
$\Lll=\diag(l_1, \ldots,l_p)$,  
are diagonal matrices with
the eigenvalues  
$\lambda_1 \ge\ldots\ge \lambda_p>0$, $l_1\ge\ldots\ge l_p>0$
of $\s$ and $\sss$, respectively. We use the
notations $\llambda=(\lambda_1,\ldots,\lambda_p)$  and $\lll=(l_1,\ldots,l_p)$ 
hereafter. By the requirement
that 
$$
\gtilde=(\widetilde{g}_{ij})=\ggammat\g
$$
has positive diagonal elements, 
the spectral decomposition $\sss=\g\Lll\ggt$ is almost surely uniquely determined. 
Then almost surely there exists a one-to-one correspondence between
the set $\{\sss \mid \sss>0\}$ and ${\cal L}\times\opp$, where
$$
{\cal L}=\{\lll \mid l_1>\cdots>l_p>0\}, \qquad
\opp=\{\gtilde \in \op \mid  \widetilde{g}_{ii}>0,\ 1\leq i \leq p\}.
$$

Let $m$ ($m_i$ in Subsection
\ref{subsec:multiblock}) denote the dividing point of the first block
and the second block of the eigenvalues.
Now we parameterize $\llambda$,$\lll$ as follows;
\begin{equation}
\label{def_tillamb}
\lambda_i=
\left\{
\begin{array}{cl}
\xi_i \aa,& \mbox{ if }i=1,\ldots, m,\\
\xi_i \bb,& \mbox{ if }i=m+1,\ldots, p,
\end{array}
\right.
\end{equation}
\begin{equation}
l_i=
\left\{
\begin{array}{cl}
\label{def_till}
d_i \aa,& \mbox{ if }i=1,\ldots, m,\\
d_i \bb,& \mbox{ if }i=m+1,\ldots, p,
\end{array}
\right.
\end{equation}
In this paper we always consider $\xi$'s are given and fixed. We also use the notations,
$$
\xxxi=\diag(\xi_1,\ldots,\xi_p),\qquad
\xxi=(\xi_1,\ldots,\xi_p),
$$
$$
\dd=\diag(d_1,\ldots,d_p),\qquad
\ddd=(d_1,\ldots,d_p).
$$
We will investigating the asymptotic distribution of $\sss$ as
$\bdiva$ goes to 0 while $\xxxi$ is fixed and its application to the estimation of $\s$. 
The state $\bdiva \approx 0$ means that the eigenvalues of $\s$ are
two-block-wise ``largely'' dispersed. In the following, the notation
$\bdiva\rightarrow 0$ means a limiting operation $n\rightarrow \infty$ with
arbitrary sequences $\alpha_n, \beta_n$, $n=1,2,\ldots$, such that
$\beta_n/\alpha_n\rightarrow 0$.

We briefly describe the content of the following sections. In
Subsection \ref{subsec:localcoord} we prepare a local coordinate
system of $\opp$ around $\ii_p$. %
In Subsection \ref{subsec:mainresult} we present our main results on
asymptotic distributions and we further discuss the case of
multi-block-wise infinite dispersion in Subsection
\ref{subsec:multiblock}.  Section \ref{sec:appli_to_est_S} deals with
the estimation of $\s$ from decision-theoretic framework.  In
Subsection \ref{subsec:intro_appli_to_est_S} we introduce orthogonally
equivariant estimators and two loss functions and in Subsection
\ref{subsec:asympt_risk} we calculate the asymptotic
risks. %
We concentrate on the special case of block-wise identity covariance
matrices in Subsection \ref{subsec:Minimum Asymptotic Risk Estimator},
which is practically important, and we propose  the best estimator
for the case with respect to each loss function.
In Subsection
\ref{subsec:numerical_analysis} the convergence speed of both
distributions and risks are numerically evaluated. Together with the
application to discriminant analysis, the numerical comparisons
show the superiority of the new estimators.  In
Appendix we present the proofs of two lemmas
and discuss
analytical calculation of the asymptotic risks.

Before concluding this subsection, we introduce some notational
conventions in this paper. In the sections other than
Subsection \ref{subsec:multiblock}, we always consider a same two-block partition
of matrices. For $\bm{A}=(a_{ij})$, a $p\times p$ matrix, 
$\bm{A}_{ij}$ $(1\leq i,j \leq 2)$ denotes the $(i,j)$-block in the 
partition
$$
\bm{A}=\left(
\begin{array}{cc}
\bm{A}_{11}&\bm{A}_{12}\\
\bm{A}_{21}&\bm{A}_{22}
\end{array}
\right),
\qquad \bm{A}_{11}:m\times m,\quad \bm{A}_{22}:(p-m)\times (p-m).
$$
If $\bm{A}$ is block diagonal, i.e.\ $\bm{A}_{12}=\bm{A}_{21}=\zz$,
we write 
$$
\bm{A}=\diag(\bm{A}_{11}, \bm{A}_{22}) = 
\left(
\begin{array}{cc}
\bm{A}_{11}& \zz\\
\zz&\bm{A}_{22}
\end{array}
\right).
$$
For the particular case of diagonal matrix
$\bm{A}=\diag(a_1,\ldots,a_p)$, we simply write $\bm{A}_1, \bm{A}_2$
instead of $\bm{A}_{11}, \bm{A}_{22}$, i.e.\ 
$\bm{A}_1=\diag(a_1,\ldots,a_m),\quad
\bm{A}_2=\diag(a_{m+1},\ldots,a_p)$.  Let $\bm{a}=(a_{ij})_{1\leq j<i
  \leq p}$ denote the vector of the elements in the lower
triangular part of $\bm{A}$, which is correspondingly partitioned as
$\bm{a}=(\bm{a}_{11}, \bm{a}_{22}, \bm{a}_{21})$, where
$$
\bm{a}_{11}=(a_{ij})_{1\leq j < i \leq m}, \quad
\bm{a}_{22}=(a_{ij})_{m+1\leq j < i\leq p}, \quad
\bm{a}_{21}=(a_{ij})_{1\leq j \leq m < i \leq p}.
$$
If $\bm{a}$ is a $p$-dimensional row vector, i.e., $\bm{a}=(a_1,\ldots,a_p)$, then we make a partition of $\bm{a}$ as
$$
\bm{a}=(\bm{a}_1,\bm{a}_2),\qquad \bm{a}_1=(a_1,\ldots,a_m),\quad \bm{a}_2=(a_{m+1},\ldots,a_p).
$$
We write $\etr X = \exp(\tr X)$ for a square matrix $X$.
\section{Asymptotic Distribution}
\label{sec:Asympto}

\subsection{Local Coordinates}
\label{subsec:localcoord}
We consider a local coordinate of $\opp$, $\uu=(u_{ij})_{1\leq j<i
  \leq p}$, around the identity matrix $\bm{I}_p$. For the proof of
the existence of such coordinate, see Appendix B of Takemura and
Sheena (2005). 
We have the following open sets $C_\epsilon, U, V$ and functions $\phi_{ij},\ 1\leq i\leq j \leq p$;
\begin{eqnarray*}
&&C_\epsilon=\{\uu \;|\; |u_{ij}|<\epsilon, 1\leq j<i \leq p \}
\subset R^{p(p-1)/2},\\
&&\bm{0}\in U\subset \bar{U} \subset C_\epsilon, \\
&&\bm{I}_p\in V \subset \opp,
\end{eqnarray*}
and $\phi_{ij}(\uu)$ is a $C^\infty$ function on $C_\epsilon$ such that $\g(\uu)=(g_{ij}(\uu))$ defined by 
\begin{equation}
\label{def_G(u)}
\left\{
\begin{array}{rcll}
g_{ij}(\uu)&=&\phi_{ij}(\uu), &  1\leq i \leq j \leq p,\\
g_{ij}(\uu)&=&u_{ij}, & 1\leq j < i \leq p,
\end{array}
\right.
\end{equation}
is a one-to-one function from $U$ onto $V$.
Using $V$ we can construct a
finite open covering of $\opp$ as follows. 
For $\hh_1\in {\cal O}^+(m),\ \hh_2\in {\cal
  O}^+(p-m)$, let 
$$
V(\hh_1,\hh_2)=\diag(\hh_1,\hh_2)V\cap \opp
=\{\g \mid  \g=\diag(\hh_1,\hh_2)\g^* \;,\;  \exists \g^*\in V\}\cap\opp.
$$
denote the open neighborhood of $\diag(\hh_1, \hh_2)$.
Let
$$
{\cal O}(m,p-m)=
\
\{ \diag(\hh_1, \hh_2) \mid \hh_1 \in {\cal O}^+(m), \ 
 \hh_2\in {\cal O}^+(p-m) \}
$$
then
$$
{\cal O}(m,p-m) \subset \bigcup_{\hh_1\in {\cal O}^+(m), \hh_2\in {\cal
    O}^+(p-m)} V(\hh_1,\hh_2) .
$$
Since ${\cal O}(m,p-m)$ is compact, we can choose a finite number of 
sets %
$
O^{(\covernum)}=V(\hh_1^{(\covernum)},\hh_2^{(\covernum)})
$,
$\covernum=1,\dots,\covertotal$,
such that $\bigcup_{\covernum=1}^\covertotal O^{(\covernum)} \supset {\cal O}(m,p-m)$.
Let $O^{(0)}=\opp \setminus {\cal O}(m,p-m)$, then we have a finite open covering $\{O^{(\covernum)}\}_{\covernum=0}^\covertotal$ of $\opp.$
We denote the partition of unity subordinate to 
$\{ O^{(\covernum)}\}_{\covernum=0}^\covertotal$ by $\{
\iota_\covernum\}_{\covernum=0}^\covertotal$. 
Namely
for each $\covernum$, $\iota_\covernum$ is a continuous function from $\opp$ to
$[0,1]$, the support of $\iota_\covernum$ is contained  in $O^{(\covernum)}$, and
$\sum_{\covernum=0}^\covertotal \iota_\covernum(G)\equiv 1$.

For $O^{(\covernum)}$, $1\leq \covernum \leq \covertotal$, we can use
$\uu$ as a local coordinate since $\g$ in $O^{(\covernum)}$ can be
uniquely expressed as $\hh^{(\covernum)}\g(\uu)$ with some $\uu$ in
$U$, where
\begin{equation}
\label{def_H(s)}
\hh^{(\covernum)}=\diag(\hh_1^{(\covernum)},\hh_2^{(\covernum)}), \qquad \covernum=1,\ldots,\covertotal.
\end{equation}
As we will see later, we do not need a local coordinate for $O^{(0)}$, since the measure of this area asymptotically vanishes.

Now we have $(\lll,\uu)$ as a local coordinate on each ${\cal L}\times
O^{(\covernum)},\ \covernum=1,\ldots,\covertotal$. We need another local coordinate to
investigate the asymptotic behavior of $\sss$. Let
$\qq=(q_{ij})_{1\leq j<i \leq p}$ be defined as follows as another coordinate on $O^{(\covernum)}$ for a fixed $\covernum$, $\covernum=1,\ldots,\covertotal$; if $1\leq j
\leq m <i \leq p,$
\begin{eqnarray}
\label{def_q_ij}
q_{ij}&=&
l_j^{1/2}\lambda_i^{-1/2}\sum_{t=m+1}^{p}(\hh_2^{(\covernum)})_{i-m,t-m}\:u_{tj}\nonumber\\
&=&\aa^{1/2}\bb^{-1/2}\:d_j^{\:1/2}\xi_i^{-1/2}\sum_{t=m+1}^{p}(\hh_2^{(\covernum)})_{i-m,t-m}
\:u_{tj}
\end{eqnarray}
and $q_{ij}=u_{ij}$ otherwise. If we use matrices $\bm{Q}=(q_{ij}),\ \bm{U}=(u_{ij})$ and their partitions, (\ref{def_q_ij}) is the same as
\begin{equation}
\label{Exp_Q}
\qqq_{21}=\aa^{1/2}\bb^{-1/2}\xxxi_2^{-1/2}\hh_2^{(\covernum)} \uuu_{21} \dd_1^{1/2},\quad
\qqq_{11}=\uuu_{11},\quad
\qqq_{22}=\uuu_{22}.
\end{equation}
Conversely
\begin{equation}
\uuu_{21}=\aa^{-1/2}\bb^{1/2}\hh_2^{(\covernum)}{}'\xxxi_2^{1/2}\qqq_{21}\dd_1^{-1/2},\quad
\uuu_{11}=\qqq_{11},\quad
\uuu_{22}=\qqq_{22},
\end{equation}
or 
\begin{equation}
\label{def_u( )}
u_{ij}=\left\{
\begin{array}{ll}
\displaystyle{\alpha^{-1/2}\beta^{1/2}\sum_{t=m+1}^{p}(\hh_2^{(\covernum)})_{t-m,i-m}\:q_{tj}\:\xi_t^{1/2}\:d_j^{-1/2},} &
\mbox{if $1\leq j \leq m <i \leq p,$}\\
q_{ij}, &\mbox{otherwise.}
\end{array}
\right.
\end{equation}
Pairing $\qq=(q_{ij})_{1\leq j<i \leq p}$ with $\ddd=(d_1,\ldots, d_p)$, we have another local coordinate $(\ddd,\qq)$ on ${\cal D}\times O^{(\covernum)}$, where
\begin{equation}
{\cal D}=({\cal D}_1\times {\cal D}_2) \cap {\cal D}_3
\end{equation}
with
$$
\begin{array}{crl}
{\cal D}_1&=&\{\ddd_1 \mid d_1 > \cdots > d_m >0\}\\
{\cal D}_2&=&\{\ddd_2 \mid d_{m+1} > \cdots > d_p >0\}\\
{\cal D}_3&=&\{(\ddd_1,\ddd_2) \mid d_m/d_{m+1}>\bdiva\}.
\end{array}
$$
The Jacobian of the transformation $J((\lll,\uu)\rightarrow (\ddd,\qq))$ is given by
\begin{eqnarray}
\label{jac_lu_tilq}
\left|
\mbox{det}\left(
\frac{\partial (\lll,\uu)}{\partial (\ddd,\qq)}
\right)
\right|&=&
\left|\mbox{det}\left(\frac{\partial \lll}{\partial \ddd}\right)\right|
\left|\mbox{det}\left(\frac{\partial \uu}{\partial \qq}\right)\right|\nonumber\\
&=&\aa^m \bb^{p-m} \prod_{j\leq  m
  <i}\Bigl(d_j^{-\frac{1}{2}}\xi_i^{\frac{1}{2}}
\aa^{-\frac{1}{2}}\bb^{\frac{1}{2}}\Bigr)\nonumber\\
&=&\aa^{m-\frac{m(p-m)}{2}}\bb^{p-m+\frac{m(p-m)}{2}}\prod_{j=1}^md_j^{-\frac{(p-m)}{2}}\prod_{i=m+1}^p\xi_i^{\frac{m}{2}}.
\end{eqnarray}
\subsection{Main Results}
\label{subsec:mainresult}
The following theorem says that $\gtilde$ asymptotically separates into two orthogonal matrices $\gtilde_{11}, \gtilde_{22}$ on the diagonal blocks.
\begin{theorem}
\label{thm:asympt_g21} \ 
\begin{itemize}
\item[{\rm 1}] As $\bdiva \rightarrow 0$, $\gtilde_{21}\plim \bm{0}$.
\item[{\rm 2}] $\lim_{\bdiva\rightarrow 0}P(\gtilde \in O)=1$
for any open set $O\subset \opp$ including ${\cal O}(m,p-m)$.
\end{itemize}
\end{theorem}
\proof
Since 2 is easily proved from 1, we only prove 1 here.  Let 
$$
\bar{\sss} =(\bar{s}_{ij})=\Llambda^{-\frac{1}{2}}\ggammat\sss\ggamma\Llambda^{-\frac{1}{2}}
=\Llambda^{-\frac{1}{2}}\gtilde\Lll\gtildet\Llambda^{-\frac{1}{2}}
\sim\ww_p(n, \iip), 
$$
Suppose $1\leq j \leq m <i\leq p$. Note that 
$$
\bar s_{ii}
  =(\widetilde g_{i1}^2l_1+\cdots+\widetilde g_{ip}^2l_p)\lambda_i^{-1}.
$$
Therefore 
\begin{equation}
\label{eq:tilde-g2-bound}
 \widetilde g_{ij}^2\le \bar s_{ii} \frac{\lambda_i}{l_j}
 =\bar s_{ii}\frac{\lambda_j}{l_j}\frac{\lambda_i}{\lambda_j}
\le \bar s_{ii}\frac{\lambda_j}{l_j}\frac{\xi_i}{\xi_j}\frac{\bb}{\aa}.
\end{equation}
Since $\bar s_{ii}$ is distributed independently of $\s$, 
for any $\epsilon>0$, there exists $M$ such that
\begin{equation}
\label{eq:tilde-wiiM}
 P(\bar s_{ii}<M)>1-\epsilon,\qquad \forall \s.
\end{equation}
Besides, from the result of Lemma 1 of Takemura \& Sheena (2005), for any $\epsilon >0$, there exists $C$ such that 
\begin{equation}
\label{eq:lemma1-consequence}
   P\left(\frac{\lambda_j}{l_j}<C\right)>1-\epsilon,\qquad\forall \s.
\end{equation}
From (\ref{eq:tilde-wiiM}) 
and (\ref{eq:lemma1-consequence}) we have
$$
\bar s_{ii}\frac{\lambda_j}{l_j}
\frac{\bb}{\aa}
\stackrel{p}{\rightarrow}0\
\quad\mbox{ as }\ \frac{\bb}{\aa}\rightarrow 0.
$$
From this fact and (\ref{eq:tilde-g2-bound}) we have
$$
  \widetilde g_{ij}^2\stackrel{p}{\rightarrow}0\quad
  \mbox{ as }\ \frac{\bb}{\aa}\rightarrow 0, 
  \qquad  1\le \forall j\leq m <\forall i\le p.
$$
\hfill\qed
\\

Next we state a rather technical lemma, which will be used in the
proofs of some theorems. Consider a random variable
$x(\g,\lll,\llambda,\aa,\bb).$ We are often interested in the
asymptotic expectation of $x(\g,\lll,\llambda,\aa,\bb)$ as
$\bdiva\rightarrow 0$ while $\ggamma$ is fixed. For fixed $\ggamma$
and $\hh^{(\covernum)}=\diag(\hh^{(\covernum)}_1,\hh^{(\covernum)}_2),\ \hh^{(\covernum)}_1\in {\cal O}^+(m),\ \hh^{(\covernum)}_2\in {\cal O}^+(p-m)$, somewhat abusing the notation, let
\begin{equation}
\label{exp_x}
x(\ddd,\qq,\xxi,\aa,\bb;\ggamma,\hh^{(\covernum)})=
x(\ggamma\hh^{(\covernum)}\g(\uu(\ddd,\qq,\xxi,\aa,\bb)),\lll(\ddd,\aa,\bb),\llambda(\xxi,\aa,\bb),\aa,\bb)
\end{equation}
for emphasizing the right-hand side as the function of $(\ddd,\qq,\xxi,\aa,\bb)$, where $\g(\uu)$, 
$\uu(\ddd,\qq,\xxi,\aa,\bb)$, $\lll(\ddd,\aa,\bb)$, $\llambda(\xxi,\aa,\bb)$ are
respectively defined by  (\ref{def_G(u)}), (\ref{def_u( )}), (\ref{def_till}) and (\ref{def_tillamb}).
For $\bm{u}=(\bm{u}_{11},\bm{u}_{22},\bm{u}_{21})$, we have
\begin{equation}
\label{asym_u}
\lim_{\bdiva \rightarrow 0}\uu(\ddd,\qq,\xxi,\aa,\bb)
=\lim_{\bdiva\rightarrow 0}(\bm{u}_{11}(\bm{q}_{11}),\bm{u}_{22}(\bm{q}_{22}),\bm{u}_{21}(\ddd,\qq,\xxi,\aa,\bb))
=(\bm{q}_{11},\bm{q}_{22},\bm{0}),
\end{equation}
hence
\begin{equation}
\label{asympt_G}
\lim_{\bdiva\rightarrow 0}\g(\bm{u}(\ddd,\qq,\xxi,\aa,\bb))=\g(\bm{q}_{11},\bm{q}_{22},\bm{0}).
\end{equation}

\begin{lem}
\label{lem:asympt_expec}
Suppose that there exist some $a<1/2$ and $b >0$ such that 
\begin{equation}
\label{bound_cond_x}
|x(\ggamma\g,\lll,\llambda,\aa,\bb)|\leq b
\etr(a\g\Lll\ggt\Llambdain)  \  \mbox{ a.e. in } (\g,\lll)
\end{equation}
and suppose that for each $\covernum,\ \covernum=1,\ldots,\covertotal$,
$
\lim_{\bdiva\rightarrow 0}x(\ddd,\qq,\xxi,\aa,\bb;\ggamma,\hh^{(\covernum)})
$
exists and equals to a function
\begin{equation}
\label{conv_x}
\bar{x}_{\ggamma}(\hh^{(\covernum)}\g(\qq_{11},\qq_{22},\bm{0}),\ddd,\qqq_{21},\xxi).
\end{equation}
Then
\begin{eqnarray}
&&\label{basi_asym_lem}
\lim_{\bdiva\rightarrow 0}E[x(\g,\lll,\llambda,\aa,\bb)]\\ \nonumber
&& \qquad=E[\bar{x}_{\ggamma}
(\diag(\g_{11}(\ww_{11}), \g_{22}(\ww_{22}))
,(\ddd_1(\ww_{11}),\ddd_2(\ww_{22})),\bm{Z}_{21},\xxi)],
\end{eqnarray}
where
the expectation on the right side of (\ref{basi_asym_lem}) is taken with respect to the following mutually independent distributions
\begin{eqnarray}
\label{asympt_dist_lemma1}
\ww_{11}&\sim &\ww_m(n,\xxxi_1),\nonumber\\
\ww_{22}&\sim &\ww_{p-m}(n-m,\xxxi_2),\\ 
\bm{Z}_{21}&\sim&\bm{N}_{(p-m)\times m}(\bm{0},\bm{I}_{p-m}\otimes \bm{I}_m),\nonumber
\end{eqnarray}
and $\g_{ss}(\ww_{ss}),\ddd_s(\ww_{ss})$, $s=1,2$, are the components
in the unique spectral decomposition of $\ww_{ss}$
for $s=1,2$;
\begin{equation}
\label{decomp_W_ii}
\begin{array}{lll}
\ww_{11}=\g_{11}\dd_1\g'_{11}, & 
\dd_1=\diag(d_1,\ldots,d_m),
 & \ddd_1=(d_1,\ldots,d_m),\\
\ww_{22}=\g_{22}\dd_2\g'_{22}, &
\dd_2=\diag(d_{m+1},\ldots,d_p),
 & \ddd_2=(d_{m+1},\ldots,d_p).
\end{array}
\end{equation}
\end{lem}

The proof is given in Appendix.  

The following theorem on the asymptotic distributions is actually a corollary of Lemma \ref{lem:asympt_expec}. Let
$$
\widetilde{\ww}_{11}=\gtilde_{11}\dd_1 \gtildet_{11},
$$
$$
\widetilde{\ww}_{22}=\gtilde_{22}\dd_2 \gtildet_{22},
$$
$$
\widetilde{\bm{Z}}_{21}=\aa^{1/2}\bb^{-1/2}\xxxi_2^{-1/2} \gtilde_{21}
\dd_1^{1/2},
$$
where all the elements on the right-hand side are defined in Section \ref{sec:intro}.
\begin{theorem}
\label{thm:asympt_dist}
As $\bdiva\rightarrow 0$,
$$
\widetilde{\ww}_{11}\dlim \ww_m(n,\xxxi_1),
$$
$$
\widetilde{\ww}_{22}\dlim \ww_{p-m}(n-m,\xxxi_2),
$$
$$
\widetilde{\bm{Z}}_{21}\dlim \bm{N}_{(p-m)\times m}(\bm{0},\bm{I}_{p-m}\otimes \bm{I}_m)
$$
and $\widetilde{\ww}_{11}$, $\widetilde{\ww}_{22}$, $\widetilde{\bm{Z}}_{21}$ are asymptotically mutually independently distributed.
\end{theorem}
\proof
Let $\bm{\Theta}_{11}: m\times m$ symmetric matrix,  
$\bm{\Theta}_{22}: (p-m)\times (p-m)$ symmetric matrix and  
$\bm{\Theta}_{21}: m \times (p-m)$ matrix.
Consider the moment generating function 
\begin{eqnarray*}
x(\g,\lll,\llambda,\aa,\bb)
&=&\exp(\tr\widetilde{\ww}_{11}\bm{\Theta}_{11}
+\tr\widetilde{\ww}_{22}\bm{\Theta}_{22}
+\tr\widetilde{\bm{Z}}_{21}\bm{\Theta}_{21}) \\
&=&
\exp(\sum_{s=1}^2 \tr\widetilde{\ww}_{ss}\bm{\Theta}_{ss}
+\tr\widetilde{\bm{Z}}_{21}\bm{\Theta}_{21}) .
\end{eqnarray*}
For $\hh^{(\covernum)}=\diag(\hh^{(\covernum)}_1,\hh^{(\covernum)}_2),\ \hh^{(\covernum)}_1\in {\cal O}^+(m),\ \hh^{(\covernum)}_2\in {\cal O}^+(p-m)$, we have
\begin{eqnarray*}
x(\ggamma\hh^{(\covernum)}\g(\uu),\lll,\llambda,\aa,\bb)
&=&\exp\Bigl\{
\sum_{s=1}^2
\tr
(\hh^{(\covernum)}\g(\uu))_{ss}\dd_s
(\hh^{(\covernum)}\g(\uu))_{ss}'
\bm{\Theta}_{ss}\\
& &\qquad +\tr \aa^{1/2}\bb^{-1/2}\xxxi_2^{-1/2}(\hh^{(\covernum)}\g(\uu))_{21}
\dd_1^{1/2}\bm{\Theta}_{21}
\Bigr\}.
\end{eqnarray*}
From (\ref{def_G(u)})%
$$
(\hh^{(\covernum)}\g(\uu))_{21}=\hh^{(\covernum)}_2\bm{U}_{21},
$$
hence from (\ref{Exp_Q})
$$
\aa^{1/2}\bb^{-1/2}\xxxi_2^{-1/2}(\hh^{(\covernum)}\g(\uu))_{21}
\dd_1^{1/2}=\bm{Q}_{21}.
$$
This leads to
$$
x(\ddd,\qq,\xxi,\aa,\bb;\ggamma,\hh^{(\covernum)})
=\exp\Bigl\{
\sum_{s=1}^2 \tr
(\hh^{(\covernum)}\g(\uu))_{ss}\dd_s(\hh^{(\covernum)}\g(\uu))_{ss}'\bm{\Theta}_{ss}\\
+\tr \qqq_{21}\bm{\Theta}_{21}
\Bigr\},
$$
with $\uu=\uu(\ddd,\qq,\xxi,\aa,\bb)$.
Therefore from (\ref{asympt_G})
\begin{eqnarray*}
\lefteqn{\lim_{\bdiva\rightarrow 0}x(\ddd,\qq,\xxi,\aa,\bb;\ggamma,\hh^{(\covernum)})}\\
&=&\exp\Bigl\{
\sum_{s=1}^2 \tr
(\hh^{(\covernum)}\g(\qq_{11},\qq_{22},\zz))_{ss}\dd_s
(\hh^{(\covernum)}\g(\qq_{11},\qq_{22},\zz))_{ss}'\bm{\Theta}_{ss}
+\tr\bm{Q}_{21}\bm{\Theta}_{21}\nonumber
\Bigr\}.
\end{eqnarray*}
From Lemma \ref{lem:asympt_expec}, 
\begin{eqnarray*}
\lefteqn{\lim_{\bdiva\rightarrow 0}
E[\exp(\tr\widetilde{\ww}_{11}\bm{\Theta}_{11}
  +\tr\widetilde{\ww}_{22}\bm{\Theta}_{22}
  +\tr\widetilde{\bm{Z}}_{21}\bm{\Theta}_{21})]}\\
&=&E[\exp\{
\sum_{s=1}^2 \tr\g_{ss}(\ww_{ss})\dd_s(\ww_{ss})\g_{ss}(\ww_{ss})'
\bm{\Theta}_{ss}
+\tr\bm{Z}_{21}\bm{\Theta}_{21}\}]\\
&=&E[\etr\ww_{11}\bm{\Theta}_{11}]E[\etr\ww_{22}\bm{\Theta}_{22}]
E[\etr\bm{Z}_{21}\bm{\Theta}_{21}],
\end{eqnarray*}
where in the second and third equations the
expectations are taken with respect to the distributions
(\ref{asympt_dist_lemma1}) in Lemma \ref{lem:asympt_expec}.  \hfill
\qed

\subsection{Multi-block Partition}
\label{subsec:multiblock}
In this section, we extend Theorem \ref{thm:asympt_dist} into multi-block cases. We partition $(1,\ldots,p)$ into $k$ blocks; 
$$
\begin{array}{rcl}
\mbox{ 1st block}&&(m_0+1,\cdots,m_1),\\
\mbox{ 2nd block}&&(m_1+1,\ldots,m_2),\\
&\vdots& \\
\mbox{ $k$th block}&&(m_{k-1}+1,\cdots,m_k),
\end{array}
$$
where
$$
m_0=0<m_1<m_2<\cdots <m_k=p.
$$
Let $[i]$, $i=1,\ldots,p$, denote the block containing $i$, i.e.,
$$
[i]=s,\quad \mbox{if } m_{s-1}+1\leq i \leq m_s.
$$
We also use the notations $\mdiff_s=m_s-m_{s-1},\ s=1,\ldots,k$, for
the block sizes.

Correspondingly to the above partition, we make the following partition of a $p\times p$ matrix $\bm{A}=(a_{ij})$;
$$
\bm{A}=
\left(
\begin{array}{ccc}
\bm{A}_{11}& \cdots & \bm{A}_{1k} \\
\vdots  & \ddots & \vdots   \\
\bm{A}_{k1}& \cdots & \bm{A}_{kk}
\end{array}
\right),\quad \bm{A}_{st}:\mdiff_s\times \mdiff_t \mbox{ matrix},\ 1\leq s,t \leq k.
$$
For a diagonal matrix $\bm{A}={\rm diag}(a_1,\ldots,a_p)$, we use the notation
$$
\bm{A}=
\left(
\begin{array}{ccc}
\bm{A}_{1} &  & \bm{0}\\
           &\ddots &   \\
\bm{0}     & &\bm{A}_{k}
\end{array}
\right),
\qquad \bm{A}_s={\rm diag}(a_{m_{s-1}+1},\ldots,a_{m_s}),\ s=1,\ldots,k.
$$
Consider the following parametrization 
of $\lll,\llambda$ 
$$
\begin{array}{rcll}
\lambda_i&=&\xi_i \aa_{[i]},\quad &1\leq i \leq p.\\
l_i&=&d_i \aa_{[i]},\quad &1\leq i \leq p,
\end{array}
$$
In this subsection we again consider that $\xi_i$'s are fixed.
Now we define $\widetilde{\ww}_{ss}$, $\widetilde{\bm{Z}}_{st}$, 
$1\leq t < s \leq k$;
\begin{eqnarray*}
\widetilde{\ww}_{ss}&=&\gtilde_{ss}\dd_s\gtildet_{ss},\\
\widetilde{\bm{Z}}_{st}&=&\alpha_t^{1/2}\alpha_s^{-1/2}\xxxi_s^{-1/2}\gtilde_{st}\dd_t^{1/2},
\end{eqnarray*}
where notations of the right-hand side are defined in Section \ref{sec:intro}.
The following theorem is the extension of Theorem \ref{thm:asympt_dist}.
\begin{theorem}
\label{thm:asympt_dist_multi}
As $(\aa_2/\aa_1, \aa_3/\aa_2,\cdots, \aa_k/\aa_{k-1})\rightarrow \bm{0}$, 
\begin{eqnarray*}
\widetilde{\ww}_{ss}&\dlim& \ww_{\alpha_s}(n-m_{s-1},\xxxi_s),\quad 1\leq s \leq k,\\
\widetilde{\bm{Z}}_{st}&\dlim& \bm{N}_{\mdiff_s\times \mdiff_t}(\bm{0},\bm{I}_{\mdiff_s}
\otimes \bm{I}_{\mdiff_t}),\quad 1\leq t < s \leq k,
\end{eqnarray*}
and $\widetilde{\ww}_{ss}(1\leq s \leq k), \widetilde{\bm{Z}}_{st}(1\leq t<s \leq k)$ are asymptotically mutually independently distributed.
\end{theorem}
\proof Though we can prove the theorem in the same manner as the proof
of Theorem \ref{thm:asympt_dist}, it is notationally too cumbersome.
Instead we will prove the theorem by using Theorem
\ref{thm:asympt_dist} recursively.  Let $r_1=\aa_1$ and
$r_t=\aa_t/\aa_{t-1},\ t=2,\ldots k,$ then $\prod_{t=1}^sr_t=\aa_s,\
s=1,\ldots,k.$ Note for $1\leq i \leq p$,
$$
l_i=d_i \aa_{[i]}=d_i\prod_{t=1}^{[i]}r_t, \qquad
\lambda_i=\xi_i \aa_{[i]}=\xi_i\prod_{t=1}^{[i]}r_t .
$$
We consider the moment generating function
$$
E\Bigl[\exp\Bigl(\tr\sum_{s=1}^k
 \widetilde{\ww}_{ss}\bm{\Theta}_{ss}
 +\tr\sum_{1\leq t<s \leq k}\widetilde{\bm{Z}}_{st}\bm{\Theta}_{st}\Bigr)\Bigr],
$$
where $\bm{\Theta}_{ss}(1\leq s \leq k)$ and $\bm{\Theta}_{st}(1\leq
t<s \leq k)$ are respectively a $\mdiff_s \times \mdiff_s$ symmetric
matrix and a 
$\mdiff_t \times \mdiff_s$ matrix. We have
\begin{eqnarray*}
\lefteqn{\lim_{(a_2/a_1,\ldots,a_k/a_{k-1})\rightarrow 0}
  E\Bigl[\exp\Bigl(\tr\sum_{s=1}^k \widetilde{\ww}_{ss}\bm{\Theta}_{ss}
  +\tr\sum_{1\leq t<s \leq k}\widetilde{\bm{Z}}_{st}\bm{\Theta}_{st}\Bigr)\Bigr]}\nonumber\\
&=&
\lim_{(r_2,\ldots,r_k)\rightarrow 0}
  E\Bigl[\exp\Bigl(\tr\sum_{s=1}^k
  \widetilde{\ww}_{ss}\bm{\Theta}_{ss}
  +\tr\sum_{1\leq t<s \leq k}\widetilde{\bm{Z}}_{st}\bm{\Theta}_{st}\Bigr)\Bigr]\nonumber\\
&=&
\lim_{r_2\rightarrow 0}\cdots \lim_{r_k\rightarrow 0}
  E\Bigl[\exp\Bigl(\tr\sum_{s=1}^k
  \widetilde{\ww}_{ss}\bm{\Theta}_{ss}
  +\tr\sum_{1\leq t<s \leq k}\widetilde{\bm{Z}}_{st}\bm{\Theta}_{st}\Bigr)\Bigr]
\end{eqnarray*}
We omit technical arguments on uniform convergences, which guarantees
the decomposition of $\lim_{(r_2,\ldots,r_k)\rightarrow 0}$ in the
second line into step by step limiting operations
$\lim_{r_2\rightarrow 0}\cdots \lim_{r_k\rightarrow 0}$ in the third
line.  

Consider the partitions;
$$
\gtilde=\left(
\begin{array}{cccl}
            &                  &                  & \gtilde_{1k}    \\
            &  \ \ \gtilde^{(k-1)} &                  & \ \vdots          \\
            &                  &                  & \gtilde_{k-1\:k}\\
\gtilde_{k1}&   \cdots         & \gtilde_{k\:k-1} & \gtilde_{kk}
\end{array}
\right)
\ \  \mbox{where}\ \ 
\gtilde^{(k-1)}=\left(
\begin{array}{ccc}
\gtilde_{11\ \ \ } & \cdots & \gtilde_{1 k-1\ \ \ }\\
\vdots & \ddots & \vdots \\
\gtilde_{k-1 1}  & \cdots & \gtilde_{k-1 k-1}
\end{array}
\right).
$$
Define $\dd^*,\ \xxxi^*$ as partitioned matrices;
$$
\dd^*=\left(
\begin{array}{cc}
\Lll^{(k-1)} & \bm{0}\\
\bm{0}\ \  & \dd_{k}\ 
\end{array}
\right),\qquad
\xxxi^*=\left(
\begin{array}{cc}
\Llambda^{(k-1)} & \bm{0}\\
\bm{0}\ \  & \xxxi_{k}\ 
\end{array}
\right),
$$
where
$$
\Lll^{(k-1)}=\mbox{diag}(l_1,\ldots,l_{m_{k-1}}),
\qquad
\Llambda^{(k-1)}=\mbox{diag}(\lambda_1,\ldots,\lambda_{m_{k-1}}).
$$
Let $\aa=1,\ \bb=\aa_k=\prod_{t=1}^k r_t.$ Then
$$
\Lll=\left(
\begin{array}{cc}
\Lll^{(k-1)}\alpha & \bm{0} \\
\bm{0}       & \dd_k \beta
\end{array}
\right),\qquad
\Llambda=\left(
\begin{array}{cc}
\Llambda^{(k-1)}\alpha & \bm{0} \\
\bm{0}           & \xxxi_k \beta
\end{array}
\right).
$$
Since as $r_k \rightarrow 0$, $\bdiva \rightarrow 0$, from Theorem \ref{thm:asympt_dist}, we have
\begin{eqnarray*}
&&\sss^{(k-1)}=\gtilde^{(k-1)} \Lll^{(k-1)} \gtilde^{(k-1)}{}'\dlim \ww_{m_{k-1}}(n,\Llambda^{(k-1)}),\\
&&\widetilde{\ww}_{kk} \dlim \ww_{\mdiff_k}(n-m_{k-1}, \xxxi_k),\\
&&\widetilde{\bm{Z}}_{kt} \dlim  \bm{N}_{\mdiff_k \times \mdiff_t}(\bm{0}, \ii_{\mdiff_k} \otimes \ii_{\mdiff_t}),\quad 1\leq t \leq k-1,
\end{eqnarray*}
and the asymptotic distributions are mutually independent. Therefore
\begin{eqnarray*}
\lefteqn{\lim_{r_k\rightarrow 0}E\Bigl[
  \exp\Bigl(\tr\sum_{s=1}^k \widetilde{\ww}_{ss}\bm{\Theta}_{ss}
  +\tr\sum_{1\leq t<s \leq k}\widetilde{\bm{Z}}_{st}\bm{\Theta}_{st}\Bigr)\Bigr]}\\
&=&E\Bigl[\exp\Bigl(\tr\sum_{s=1}^{k-1}
\widetilde{\ww}_{ss}(\sss^{(k-1)})\bm{\Theta}_{ss}
+\tr \sum_{1\leq t<s \leq k-1}\widetilde{\bm{Z}}_{st}(\sss^{(k-1)})\bm{\Theta}_{st}\Bigr)\Bigr]\\
& &\qquad \times E\Bigl[\etr\widetilde{\ww}_{kk}\bm{\Theta}_{kk}\Bigr]
\times \prod_{t=1}^{k-1}E\Bigl[\etr\widetilde{\bm{Z}}_{kt}\bm{\Theta}_{kt}\Bigr],
\end{eqnarray*}
where the expectations on the right-hand side is taken with respect to the above asymptotic distributions. If we apply Theorem \ref{thm:asympt_dist} again to $\sss^{(k-1)}$ and recursively to the upper-left block Wishart distribution which asymptotically arises, we gain the result. \hfill \qed
\\

Note that Theorem \ref{thm:asympt_dist_multi} reduces 
to Theorem 2 of Takemura and Sheena (2005) for the extreme case 
of 1-element blocks $\mdiff_s=1$, $s=1,\ldots,p$.  Therefore
Theorem \ref{thm:asympt_dist_multi} is a generalization of Theorem 2 of Takemura and Sheena (2005).
\section{Application to Estimation of $\s$}
\label{sec:appli_to_est_S}
\subsection{Loss Functions and Orthogonally Equivariant Estimators}
\label{subsec:intro_appli_to_est_S}
In this section, we apply the asymptotic result on the distribution of $\sss$ to the estimation of $\s$ 
when $\bdiva$ vanishes. We take a decision-theoretic approach to evaluate the performance of the estimators. 
We deal with the two loss functions; one is Stein's loss (entropy loss) function
\begin{equation}
\label{def_stein_loss}
L_1(\sh,\s)=\tr(\sh\si)-\log|\sh\si|-p,
\end{equation}
and the other is a scale-invariant quadratic loss function
\begin{equation}
\label{def_quad_loss}
L_2(\sh,\s)=\tr(\sh\si-\ii_p)^2.
\end{equation}
The associated risk functions are denoted as
$$
R_\lossd(\sh,\s)=E[L_\lossd(\sh,\s)], \quad \lossd=1,2.
$$

The classical estimator of $\s$ is the unbiased estimator %
$$
\sh^U=n^{-1}\sss,
$$
which has been widely used for many statistical analysis, especially
with statistical software packages.  However, as James and Stein
(1961) showed, this estimator is neither minimax nor admissible with
Stein's loss function (\ref{def_stein_loss}). The same drawback with
respect to the quadratic loss function (\ref{def_quad_loss}) was
reported by Olkin and Selliah (1977). Following these initiative
papers, much literature has been written seeking for a superior
estimator to $\sh^U$. 
See Pal (1993) for the review on the
estimation of $\s$. 
In this paper we only refer to orthogonally equivariant 
estimators proposed by Stein (1977), Dey and Srinivasan
(1985) and Krishnamoorthy and Gupta (1989).  An estimator of the form
$$
\sh=\g\mbox{\boldmath $\Psi$}(\Lll)\ggt,\qquad\mbox{\boldmath $\Psi$}(\Lll)=
\diag(\psi_1(\lll),\ldots,\psi_p(\lll))
$$
is called {\it orthogonally equivariant}; i.e., 
$\sh(\g \sss \g')=\g \sh(\sss)\g', \ \forall \g\in \op$.

Stein (1977) and Dey and Srinivasan (1985) proposed the orthogonally equivariant estimator, $\sh^{SDS}$, defined by
$$
\psi_i(\lll)=l_i\Delta_i^{JS},\qquad
 1\le i \le p,
$$
where $\Delta_i^{JS}=(n+p+1-2i)^{-1}$. $\sh^{SDS}$ is of simple form
but dominates $\sh^U$ with substantially better risk w.r.t\ 
the loss function (\ref{def_stein_loss}). It is also a minimax
estimator. See Dey and Srinivasan (1985) and Sugiura and Ishibayashi
(1997) for more details. Order preservation among $\psi_i(\lll), \
i=1,\ldots,p,$ is discussed in Sheena and Takemura (1992).

The orthogonally equivariant estimator $\sh^{KG}$ is defined by
$$
\psi_i(\lll)=l_i\Delta_i^{OS},\qquad
 1\le i \le p,
$$
where $\Delta_i^{OS}$ is given by
$$
(\Delta_1^{OS},\ldots,\Delta_p^{OS})'
= \bm{A}^{-1}\bm{b}
$$
with a $p \times p$ matrix $\mbox{\boldmath $A$}=(a_{ij})$  and
a $p\times 1$ vector $\mbox{\boldmath $b$}=(b_i)$ defined by
\begin{eqnarray*}
a_{ij}&=&\left\{
\begin{array}{ll}
(n+p-2i+1)(n+p-2i+3),&\mbox{ if\quad  $i=j$,}\\
(n+p-2i+1),          &\mbox{ if\quad  $i>j$,}\\
(n+p-2j+1),          &\mbox{ if\quad  $j>i$,}
\end{array}
\right.
\\
b_i&=&n+p+1-2i,\quad i=1,\ldots,p.
\end{eqnarray*}
$\sh^{KG}$ is conjectured to be a minimax estimator which dominates $\sh^U$ w.r.t.~the loss function (\ref{def_quad_loss}). This was proved by Sheena (2002) for the case $p=2$.

In this section we only consider orthogonally equivariant estimators given by 

\begin{equation}
\label{family_orthog_equiv_est}
\psi_i(\lll)=c_il_i,\qquad
 1\le i \le p
\end{equation}
with some constant $c_i\ (1\leq i \leq p)$, or in the matrix expression,
$$
\bm{\Psi}(\Lll)=\Lllh \cc \Lllh,\qquad \cc=\diag(c_1,\ldots,c_p).
$$
It is interesting that $\sh^{SDS}$ and $\sh^{KG}$ are also the minimum risk estimators among the estimators of the form (\ref{family_orthog_equiv_est}) respectively for $L_1(\cdot,\cdot)$ and $L_2(\cdot,\cdot)$ when all the population eigenvalues are dispersed. See Takemura and Sheena (2005) for more details.
\subsection{Asymptotic Risk}
\label{subsec:asympt_risk}
This subsection is devoted to the calculation of the asymptotic risks $\widetilde{R}_\lossd(\sh,\s)$ 
$$
\widetilde{R}_\lossd(\sh,\s)=\lim_{\bdiva\rightarrow 0}R_\lossd(\sh,\s),
\qquad \lossd=1,2,
$$
for an orthogonally equivariant estimator defined by $(\ref{family_orthog_equiv_est}).$ 
Note that
\begin{eqnarray}
\label{loss_1}
R_1(\sh,\s)
&=&E[\tr \g\Lllh\cc\Lllh\ggt\ggamma\Llambdain\ggammat]
-\log|\cc|-E[\log|\s^{-1/2}\sss\s^{-1/2}|]-p\nonumber\\
&=&E[\tr\g\Lllh\cc\Lllh\ggt\ggamma\Llambdain\ggammat]
-\sum_{i=1}^p\log c_i
-\sum_{i=1}^pE[\log\chi^2_{n-i+1}]-p.\\
\label{loss_2}
R_2(\sh,\s)
&=&E[\tr(\g\Lllh\cc\Lllh\ggt\ggamma\Llambdain\ggammat-\ii_p)^2]\nonumber\\
&=&E[\tr(\g\Lllh\cc\Lllh\ggt\ggamma\Llambdain\ggammat)^2]
-2E[\tr\g\Lllh\cc\Lllh\ggt\ggamma\Llambdain\ggammat]+p.
\end{eqnarray}
For the evaluation $E[\log|\s^{-1/2}\sss\s^{-1/2}|]$, see e.g.\ (10) in p.132 of Muirhead (1982).

We start with the following lemma, the proof of which is given in Appendix.
\begin{lem}
\begin{eqnarray}
\label{asym_exp_tr_GLCLG'S}
\lefteqn{\lim_{\bdiva \rightarrow 0}
E[\tr\g\Lllh\cc\Lllh\ggt\ggamma\Llambdain\ggammat]}\nonumber\\
&=&E[\tr \g_{11} \ddh_1 \cc_1\ddh_1 \ggt_{11} 
\xxxi_1^{-1}]
+E[\tr\g_{22} \ddh_2 \cc_2\ddh_2 \ggt_{22} \xxxi_2^{-1}]\nonumber\\
& &\ +(p-m)\tr\cc_1,
\end{eqnarray}
\begin{eqnarray}
\label{asym_exp_tr^2_GLCLG'S}
\lefteqn{\lim_{\bdiva \rightarrow 0}E[\tr(\g\Lllh\cc\Lllh\ggt\ggamma\Llambdain\ggammat)^2]}\nonumber\\
&=&E[\tr(\g_{11} \ddh_1 \cc_1\ddh_1 \ggt_{11} \xxxi_1^{-1})^2]
+E[\tr(\g_{22} \ddh_2 \cc_2\ddh_2 \ggt_{22} \xxxi_2^{-1})^2]\nonumber\\
& &+2(p-m)E[\tr\cc_1^2\ddh_1 \ggt_{11} \xxxi_1^{-1}\g_{11} \ddh_1 ]+2\tr\cc_1E[\tr\xxxi_2^{-1}\g_{22} \ddh_2 \cc_2
\ddh_2 \ggt_{22} ]\nonumber\\
& &+(p-m)(p-m+2)\sum_{i=1}^m c_i^2+2(p-m)\sum_{1\leq i<s \leq m}c_ic_s,
\end{eqnarray}
where the expectations on the right-hand side in (\ref{asym_exp_tr_GLCLG'S}) and (\ref{asym_exp_tr^2_GLCLG'S}) are taken with respect to the distributions in (\ref{asympt_dist_lemma1}) and the decompositions in (\ref{decomp_W_ii}).
\end{lem}

Now suppose that under the distribution of $\ww_{ss}$, $s=1,2$, in
(\ref{asympt_dist_lemma1}) and their spectral decomposition in
(\ref{decomp_W_ii}), we estimate $\xxxi_s$, $s=1,2$, by
the following orthogonally equivariant estimators
\begin{eqnarray*}
\xxxih_1&=&\g_{11}\ddh_1\cc_1\ddh_1\ggt_{11},\qquad \cc_1=\diag(c_1,\ldots,c_m),\\
\xxxih_2&=&\g_{22}\ddh_2\cc_2\ddh_2\ggt_{22},\qquad \cc_2=\diag(c_{m+1},\ldots,c_p),
\end{eqnarray*}
then the risks w.r.t.\  each loss function (\ref{def_stein_loss}), (\ref{def_quad_loss}) are given by
\begin{eqnarray*}
R_{11}(\xxxih_1,\xxxi_1)&=&%
                  E[\tr(\xxxih_1\xxxii_1)-\log|\xxxih_1\xxxii_1|-m],\\
R_{21}(\xxxih_2,\xxxi_2)&=&%
                  E[\tr(\xxxih_2\xxxii_2)-\log|\xxxih_2\xxxii_2|-p+m],\\
R_{12}(\xxxih_1,\xxxi_1)&=&%
                  E[\tr(\xxxih_1\xxxii_1-\ii_m)^2],\\
R_{22}(\xxxih_2,\xxxi_2)&=&%
                  E[\tr(\xxxih_2\xxxii_2-\ii_{p-m})^2].
\end{eqnarray*}
The following theorem gives the decomposition of the asymptotic risk,
$\widetilde{R}_\lossd(\sh,\s)$, into the risks $R_{1\lossd},R_{2\lossd}$ and the
residuals $R_{3\lossd}$ for $\lossd=1,2.$
\begin{theorem} 
\label{theorem:asymto_risk}
For $\lossd=1,2$,
$$
\widetilde{R}_\lossd(\sh,\s)
=R_{1\lossd}(\xxxih_1,\xxxi_1)+R_{2\lossd}(\xxxih_2,\xxxi_2)
+R_{3\lossd},
$$
where 
$$
R_{31}=(p-m)\sum_{i=1}^m c_i,
$$
and
\begin{eqnarray*}
R_{32}&=&2(p-m)E[\tr\cc_1^2\ddh_1\ggt_{11}\xxxi_1^{-1}\g_{11}\ddh_1]+2\tr\cc_1E[\tr\xxxi_2^{-1}\g_{22}\ddh_2\cc_2
\ddh_2\ggt_{22}]\nonumber\\
& &+(p-m)(p-m+2)\sum_{i=1}^m c_i^2+2(p-m)\sum_{1\leq i<s \leq m}c_ic_s
-2(p-m)\sum_{i=1}^m c_i.
\end{eqnarray*}
All the expectations are taken with respect to the distributions (\ref{asympt_dist_lemma1}) and the decompositions (\ref{decomp_W_ii}).
\end{theorem}
\proof
From (\ref{loss_1}), 
\begin{eqnarray*}
R_{11}(\sh_1,\s_1)
&=&E[\tr\g_{11}\ddh_1\cc_1\ddh_1\ggt_{11}\xxxi_1^{-1}]-\sum_{i=1}^m\log{c_i}-
\sum_{i=1}^mE[\log\chi^2_{n-i+1}]-m,\\
R_{21}(\sh_2,\s_2)
&=&E[\tr\g_{22}\ddh_2\cc_2\ddh_2\ggt_{22}\xxxi_2^{-1}]-\sum_{i=m+1}^p\log{c_i}
 -\sum_{i=m+1}^{p}E[\log\chi^2_{n-i+1}]-p+m.
\end{eqnarray*}
Using (\ref{asym_exp_tr_GLCLG'S}) together with the above result, we have the result for $\widetilde{R}_1(\sh,\s).$
From (\ref{loss_2}),
\begin{eqnarray*}
R_{12}(\sh_1,\s_1)
&=&E[\tr(\g_{11}\ddh_1\cc_1\ddh_1\ggt_{11}\xxxi_1^{-1})^2]
-2E[\tr\g_{11}\ddh_1\cc_1\ddh_1\ggt_{11}\xxxi_1^{-1}]+m,\\
R_{22}(\sh_2,\s_2)
&=&E[\tr(\g_{22}\ddh_2\cc_2\ddh_2\ggt_{22}\xxxi_2^{-1})^2]
-2E[\tr\g_{22}\ddh_2\cc_2\ddh_2\ggt_{22}\xxxi_2^{-1}]+p-m.
\end{eqnarray*}
Using (\ref{asym_exp_tr^2_GLCLG'S}) and (\ref{asym_exp_tr_GLCLG'S}) together with the above result, we have the result for $\widetilde{R}_2(\sh,\s).$\hfill \qed
\subsection{Minimum Asymptotic Risk Estimator}
\label{subsec:Minimum Asymptotic Risk Estimator}
Consider the model (\ref{basic_model}) and suppose $\tau_1=\cdots
=\tau_m(=\tau)$ in (\ref{structure_S}). Then $\aa=\tau+\sigma^2$ and
$\bb=\sigma^2$ and
\begin{equation}
\label{pract_case}
\xxxi_1=\ii_m,\qquad \xxxi_2=\ii_{p-m}.
\end{equation}
This assumption may not be very realistic.  However
note that it is trivially satisfied in the one-factor model $m=1$,
which is frequently used in practice.
In this subsection we focus on the
estimation of $\s$ under the condition (\ref{pract_case}). In this
case, since we have no unknown parameters anymore, the asymptotic risk
is uniquely determined, hence we can derive the ``best'' i.e., minimum
asymptotic risk estimator among the orthogonally equivariant
estimators of the form (\ref{family_orthog_equiv_est}). The following
theorem gives the asymptotic risk for the case (\ref{pract_case}).
\begin{theorem}
\label{asym_risk_I&I}
If $\xxxi_1=\ii_m,\ \xxxi_2=\ii_{p-m}$, then the asymptotic risk
$\widetilde{R}_\lossd(\hat{\s},\s)$, $\lossd=1,2$, is given by the following
function of $\bm{c}=(c_1,\ldots,c_p)'$.
\begin{eqnarray}
\label{asympt_risk_special_1}
\widetilde{R}_1(\sh,\s)
&=&\sumi(b_ic_i-\log c_i)-\sumi E[\log \chi^2_{n-i+1}]-p,\\
\label{asympt_risk_special_2}
\widetilde{R}_2(\sh,\s)
&=&\bm{c}'\bm{A}\bm{c}-2\bm{b}'\bm{c}+p,
\end{eqnarray}
where $\bm{b}=(b_1,\ldots,b_p)'$ is given by
$$
b_i=
\left\{
\begin{array}{ll}
E[d_i]+p-m, &\mbox{\rm if}\ 1\leq i \leq m,\\
E[d_i], &\mbox{\rm if}\ m+1\leq i \leq p,
\end{array}
\right.
$$
and $p\times p$ symmetric matrix $\bm{A}=(a_{ij})$ is given by
$$
a_{ij}=
\left\{
\begin{array}{ll}
E[d_i^2+2(p-m)d_i]+(p-m)(p-m+2),&\mbox{\rm if }1\leq i=j \leq m,\\
E[d_i^2],&\mbox{\rm if }m+1 \leq i=j \leq p,\\
p-m,&\mbox{\rm if }1\leq i\ne j \leq m,\\
E[d_j],&\mbox{\rm if }1\leq i\leq m <j \leq p,\\
E[d_i],&\mbox{\rm if }1\leq j\leq m <i \leq p,\\
0,&\mbox{\rm otherwise}.
\end{array}
\right.
$$
All the expectations are taken with respect to the distribution (\ref{asympt_dist_lemma1}) and the decompositions (\ref{decomp_W_ii}) with $\xxxi_1=\bm{I}_m,\ \xxxi_2=\bm{I}_{p-m}$.
\end{theorem}
\proof
Evaluating $R_{j\lossd}(\sh_j,\s_j),\ 1\leq j,\lossd \leq 2$ in Theorem \ref{theorem:asymto_risk} when 
$\xxxi_1=\ii_m,\ \xxxi_2=\ii_{p-m},$ we have the following results.
\begin{eqnarray*}
R_{11}(\xxxih_1,\xxxi_1)&=&E[L_1(\xxxih_1,\ii_m)]
=E[\tr \xxxih_1-\log |\xxxih_1|-m]\nonumber\\
&=&E\Bigl[\sum_{i=1}^md_ic_i-\log|\ww_{11}|\Bigr]-\sum_{i=1}^m\log c_i-m\nonumber\\
&=&\sum_{i=1}^m E[d_i]c_i-E[\log |\ww_{11}|]-\sum_{i=1}^m \log c_i-m.
\nonumber\\
R_{21}(\xxxih_2,\xxxi_2)%
&=&\sum_{i=m+1}^pE[d_i]c_i-E[\log|\ww_{22}|]-\sum_{i=m+1}^p\log
c_i-p+m.
\nonumber\\
R_{12}(\xxxih_1,\xxxi_1)&=&E[L_2(\xxxih_1,\ii_m)]
=E[\tr(\xxxih_1-\ii_m)^2]\nonumber\\
&=&E[\tr\xxxih_1^2-2\tr \xxxih_1]+m
=E\Bigl[\sum_{i=1}^md_i^{\:2}c_i^2-2\sum_{i=1}^md_ic_i\Bigr]+m\nonumber\\
&=&\sum_{i=1}^mE[d_i^{\:2}]c_i^2-2\sum_{i=1}^mE[d_i]c_i+m.
\nonumber \\
R_{22}(\xxxih_2,\xxxi_2)%
&=&\sum_{i=m+1}^pE[d_i^{\:2}]c_i^2-2\sum_{i=m+1}^pE[d_i]c_i+p-m.
\nonumber 
\end{eqnarray*}
Next we calculate $R_{32}$ in Theorem \ref{theorem:asymto_risk} when $\xxxi_1=\ii_m,$ $\xxxi_2=\ii_{p-m}.$
Note that
\begin{eqnarray*}
2(p-m)E[\tr\cc_1^2\dd_1^{1/2}\g_{11}'\xxxi_1^{-1}\g_{11}\dd_1^{1/2}]&=
&2(p-m)E[\tr\cc_1^2\dd_1]\\
&=&2(p-m)\sum_{i=1}^mE[d_i]c_i^2,\\
2\tr\cc_1E[\tr\xxxi_2^{-1}\g_{22}\dd_2^{1/2}\cc_2\dd_2^{1/2}\g_{22}']
&=&2\Bigl(\sum_{i=1}^mc_i\Bigr)\Bigl(\sum_{i=m+1}^pE[d_i]c_i\Bigr).
\end{eqnarray*}
Therefore
\begin{eqnarray*}
R_{32}&=&\sum_{i=1}^mc_i^2\{(p-m)(p-m+2)+2(p-m)E[d_i]\}\\
& &\quad+2(p-m)\sum_{1\leq i<s \leq m}c_ic_s+2\sum_{1\leq i \leq m < s \leq p}c_ic_sE[d_s]-2(p-m)\sum_{i=1}^mc_i.
\nonumber
\end{eqnarray*}
Combining above results, we see that (\ref{asympt_risk_special_1}) and (\ref{asympt_risk_special_2}) hold.
\hfill \qed
\begin{cor}
The minimum asymptotic risk with respect to the loss function $L_1(\cdot,\cdot)$ is given by 
$$
\sumi \log b_i-\sumi E[\log \chi^2_{n-i-1}].
$$
It is attained by $\sh^{MA_1}$ given by $c_i=b_i^{-1},\ i=1,\dots,p$.  
The minimum asymptotic risk with respect to the loss function $L_2(\cdot,\cdot)$ is given by
$$
p-\bm{b}'\bm{A}^{-1}\bm{b}.
$$
It is attained by $\sh^{MA_2}$ given by $\bm{c}=\bm{A}^{-1}\bm{b}.$
\end{cor}
\proof
The results are easily obtained  by the minimization
$\sumi (b_i c_i -\log c_i)$ or $\bm{c}'\bm{A}\bm{c}-2\bm{b}'\bm{c}.$
\hfill\qed

\bigskip
The calculation of the asymptotic risks in Theorem \ref{asym_risk_I&I}
and the $c_i$'s of $\sh^{MA_1}$ and $\sh^{MA_2}$ requires the
evaluation of $E[d_i],  E[d_i^{\:2}], i=1,\ldots,p$, that is, the
first and the second moment of the eigenvalues of the Wishart
distribution with the identity covariance matrix. Generally we need to make
use of Monte Carlo simulation or numerical integration for the evaluation of the moments of the eigenvalues. However when $p$ is small and $n$ is appropriately even or odd depending on $p$, the analytic evaluation is feasible. See Section \ref{susec:Anal_Eval_Asymp_Risk} in Appendix for this evaluation. 

\begin{table}[btp]
\caption{$p=3, m=1$} \label{p=3, m=1}
\begin{center}
\scriptsize
\begin{tabular}{|l|l|l|l|l|l|} \hline
$n=4$ & $\sh^{U}$ & $\sh^{SDS}$ & $\sh^{KG}$ & $\sh^{MA_1}$ & $\sh^{MA_2}$ \\ \hline
$c_1$ & 0.2500 & 0.1667 & 0.1060 & 0.1667 & 0.1019 \\ \hline
$c_2$ & 0.2500 & 0.2500 & 0.1332 & 0.2000 & 0.1321 \\ \hline
$c_3$ & 0.2500 & 0.5000 & 0.1902 & 1.0000 & 0.4491 \\ \hline
Asy.Risk1 & 2.1969 & 1.6592 &  & 1.4392 &  \\ \hline
R.R.R. &  & 24.47 &  & 34.49 &  \\ \hline
Asy.Risk2 & 3.0000 &  & 1.4120 &  & 1.2792 \\ \hline
R.R.R. &  &  & 52.93 &  & 57.36 \\ \hline
\end{tabular}
\begin{tabular}{|l|l|l|l|l|l|} \hline
$n=6$ & $\sh^{U}$ & $\sh^{SDS}$ & $\sh^{KG}$ & $\sh^{MA_1}$ & $\sh^{MA_2}$ \\ \hline
$c_1$ & 0.1667 & 0.1250 & 0.0855 & 0.1250 & 0.0828 \\ \hline
$c_2$ & 0.1667 & 0.1667 & 0.1030 & 0.1304 & 0.0977 \\ \hline
$c_3$ & 0.1667 & 0.2500 & 0.1352 & 0.4286 & 0.2675 \\ \hline
Asy.Risk1 & 1.2387 & 0.9820 &  & 0.8270 &  \\ \hline
R.R.R. &  & 20.72 &  & 33.23 &  \\ \hline
Asy.Risk2 & 2.0000 &  & 1.1056 &  & 0.9644 \\ \hline
R.R.R. &  &  & 44.72 &  & 51.78 \\ \hline
\end{tabular}

\medskip

\begin{tabular}{|l|l|l|l|l|l|} \hline
$n=8$ & $\sh^{U}$ & $\sh^{SDS}$ & $\sh^{KG}$ & $\sh^{MA_1}$ & $\sh^{MA_2}$ \\ \hline
$c_1$ & 0.1250 & 0.1000 & 0.0724 & 0.1000 & 0.0707 \\ \hline
$c_2$ & 0.1250 & 0.1250 & 0.0849 & 0.0980 & 0.0782 \\ \hline
$c_3$ & 0.1250 & 0.1667 & 0.1053 & 0.2632 & 0.1878 \\ \hline
Asy.Risk1 & 0.8749 & 0.7187 &  & 0.5966 &  \\ \hline
R.R.R. &  & 17.85 &  & 31.81 &  \\ \hline
Asy.Risk2 & 1.5000 &  & 0.9140 &  & 0.7812 \\ \hline
R.R.R. &  &  & 39.07 &  & 47.92 \\ \hline
\end{tabular}
\begin{tabular}{|l|l|l|l|l|l|} \hline
$n=10$ & $\sh^{U}$ & $\sh^{SDS}$ & $\sh^{KG}$ & $\sh^{MA_1}$ & $\sh^{MA_2}$ \\ \hline
$c_1$ & 0.1000 & 0.0833 & 0.0630 & 0.0833 & 0.0619 \\ \hline
$c_2$ & 0.1000 & 0.1000 & 0.0723 & 0.0790 & 0.0655 \\ \hline
$c_3$ & 0.1000 & 0.1250 & 0.0865 & 0.1872 & 0.1438 \\ \hline
Asy.Risk1 & 0.6765 & 0.5692 &  & 0.4676 &  \\ \hline
R.R.R. &  & 15.85 &  & 30.88 &  \\ \hline
Asy.Risk2 & 1.2000 &  & 0.7817 &  & 0.6591 \\ \hline
R.R.R. &  &  & 34.86 &  & 45.07 \\ \hline
\end{tabular}

\medskip

\begin{tabular}{|l|l|l|l|l|l|} \hline
$n=20$ & $\sh^{U}$ & $\sh^{SDS}$ & $\sh^{KG}$ & $\sh^{MA_1}$ & $\sh^{MA_2}$ \\ \hline
$c_1$ & 0.0500 & 0.0455 & 0.0385 & 0.0455 & 0.0383 \\ \hline
$c_2$ & 0.0500 & 0.0500 & 0.0418 & 0.0410 & 0.0370 \\ \hline
$c_3$ & 0.0500 & 0.0556 & 0.0460 & 0.0735 & 0.0647 \\ \hline
Asy.Risk1 & 0.3164 & 0.2819 &  & 0.2251 &  \\ \hline
R.R.R. &  & 10.89 &  & 28.84 &  \\ \hline
Asy.Risk2 & 0.6000 &  & 0.4598 &  & 0.3745 \\ \hline
R.R.R. &  &  & 23.37 &  & 37.59 \\ \hline
\end{tabular}
\begin{tabular}{|l|l|l|l|l|l|} \hline
$n=50$ & $\sh^{U}$ & $\sh^{SDS}$ & $\sh^{KG}$ & $\sh^{MA_1}$ & $\sh^{MA_2}$ \\ \hline
$c_1$ & 0.0200 & 0.0192 & 0.0179 & 0.0192 & 0.0178 \\ \hline
$c_2$ & 0.0200 & 0.0200 & 0.0185 & 0.0173 & 0.0166 \\ \hline
$c_3$ & 0.0200 & 0.0208 & 0.0193 & 0.0248 & 0.0236 \\ \hline
Asy.Risk1 & 0.1236 & 0.1155 &  & 0.0901 &  \\ \hline
R.R.R. &  & 6.51 &  & 27.05 &  \\ \hline
Asy.Risk2 & 0.2400 &  & 0.2093 &  & 0.1647 \\ \hline
R.R.R. &  &  & 12.79 &  & 31.39 \\ \hline
\end{tabular}
\end{center}
\end{table} %
\begin{table}[tbhp]
\caption{$p=3, m=2$} \label{p=3, m=2}
\begin{center}
\scriptsize
\begin{tabular}{|l|l|l|l|l|l|} \hline
$n=5$ & $\sh^{U}$ & $\sh^{SDS}$ & $\sh^{KG}$ & $\sh^{MA_1}$ & $\sh^{MA_2}$ \\ \hline
$c_1$ & 0.2000 & 0.1429 & 0.0944 & 0.1154 & 0.0891 \\ \hline
$c_2$ & 0.2000 & 0.2000 & 0.1158 & 0.3000 & 0.1791 \\ \hline
$c_3$ & 0.2000 & 0.3333 & 0.1580 & 0.3333 & 0.1464 \\ \hline
Asy.Risk1 & 1.5769 & 1.3073 &  & 1.2107 &  \\ \hline
R.R.R. &  & 17.10 &  & 23.23 &  \\ \hline
Asy.Risk2 & 2.4000 &  & 1.2543 &  & 1.1919 \\ \hline
R.R.R. &  &  & 47.74 &  & 50.34 \\ \hline
\end{tabular}
\begin{tabular}{|l|l|l|l|l|l|} \hline
$n=7$ & $\sh^{U}$ & $\sh^{SDS}$ & $\sh^{KG}$ & $\sh^{MA_1}$ & $\sh^{MA_2}$ \\ \hline
$c_1$ & 0.1429 & 0.1111 & 0.0784 & 0.0893 & 0.0726 \\ \hline
$c_2$ & 0.1429 & 0.1429 & 0.0930 & 0.2083 & 0.1417 \\ \hline
$c_3$ & 0.1429 & 0.2000 & 0.1184 & 0.2000 & 0.1122 \\ \hline
Asy.Risk1 & 1.0238 & 0.8688 &  & 0.7801 &  \\ \hline
R.R.R. &  & 15.14 &  & 23.81 &  \\ \hline
Asy.Risk2 & 1.7143 &  & 1.0182 &  & 0.9455 \\ \hline
R.R.R. &  &  & 40.61 &  & 44.84 \\ \hline
\end{tabular}

\medskip

\begin{tabular}{|l|l|l|l|l|l|} \hline
$n=9$ & $\sh^{U}$ & $\sh^{SDS}$ & $\sh^{KG}$ & $\sh^{MA_1}$ & $\sh^{MA_2}$ \\ \hline
$c_1$ & 0.1111 & 0.0909 & 0.0674 & 0.0732 & 0.0616 \\ \hline
$c_2$ & 0.1111 & 0.1111 & 0.0781 & 0.1577 & 0.1162 \\ \hline
$c_3$ & 0.1111 & 0.1429 & 0.0950 & 0.1429 & 0.0914 \\ \hline
Asy.Risk1 & 0.7635 & 0.6592 &  & 0.5793 &  \\ \hline
R.R.R. &  & 13.66 &  & 24.12 &  \\ \hline
Asy.Risk2 & 1.3333 &  & 0.8585 &  & 0.7821 \\ \hline
R.R.R. &  &  & 35.61 &  & 41.34 \\ \hline
\end{tabular}
\begin{tabular}{|l|l|l|l|l|l|} \hline
$n=11$ & $\sh^{U}$ & $\sh^{SDS}$ & $\sh^{KG}$ & $\sh^{MA_1}$ & $\sh^{MA_2}$\\ \hline
$c_1$ & 0.0909 & 0.0769 & 0.0592 & 0.0623 & 0.0537\\ \hline
$c_2$ & 0.0909 & 0.0909 & 0.0674 & 0.1260 & 0.0980\\ \hline
$c_3$ & 0.0909 & 0.1111 & 0.0794 & 0.1111 & 0.0771\\ \hline
Asy.Risk1 & 0.6107 & 0.5342 &  & 0.4622 & \\ \hline
R.R.R. &  & 12.52 &  & 24.31 & \\ \hline
Asy.Risk2 & 1.0909 &  & 0.7430 &  & 0.6663\\ \hline
R.R.R. &  &  & 31.89 &  & 38.93\\ \hline
\end{tabular}

\medskip

\begin{tabular}{|l|l|l|l|l|l|} \hline
$n=21$ & $\sh^{U}$ & $\sh^{SDS}$ & $\sh^{KG}$ & $\sh^{MA_1}$ & $\sh^{MA_2}$\\ \hline
$c_1$ & 0.0476 & 0.0435 & 0.0371 & 0.0361 & 0.0331\\ \hline
$c_2$ & 0.0476 & 0.0476 & 0.0401 & 0.0613 & 0.0537\\ \hline
$c_3$ & 0.0476 & 0.0526 & 0.0439 & 0.0526 & 0.0435\\ \hline
Asy.Risk1 & 0.3040 & 0.2755 &  & 0.2281 & \\ \hline
R.R.R. &  & 9.37 &  & 24.97 & \\ \hline
Asy.Risk2 & 0.5714 &  & 0.4473 &  & 0.3815\\ \hline
R.R.R. &  &  & 21.71 &  & 33.24\\ \hline
\end{tabular}
\begin{tabular}{|l|l|l|l|l|l|} \hline
$n=51$ & $\sh^{U}$ & $\sh^{SDS}$ & $\sh^{KG}$ & $\sh^{MA_1}$ & $\sh^{MA_2}$\\ \hline
$c_1$ & 0.0196 & 0.0189 & 0.0175 & 0.0164 & 0.0158\\ \hline
$c_2$ & 0.0196 & 0.0196 & 0.0182 & 0.0232 & 0.0220\\ \hline
$c_3$ & 0.0196 & 0.0204 & 0.0189 & 0.0204 & 0.0189\\ \hline
Asy.Risk1 & 0.1191 & 0.1117 &  & 0.0881 & \\ \hline
R.R.R. &  & 6.21 &  & 25.99 & \\ \hline
Asy.Risk2 & 0.2353 &  & 0.2066 &  & 0.1666\\ \hline
R.R.R. &  &  & 12.20 &  & 29.18\\ \hline
\end{tabular}
\end{center}
\end{table} %
\begin{table}[htbp]
\caption{$p=4, m=1$} \label{p=4, m=1}
\begin{center}
\scriptsize
\begin{tabular}{|l|l|l|l|l|l|} \hline
$n=5$ & $\sh^{U}$ & $\sh^{SDS}$ & $\sh^{KG}$ & $\sh^{MA_1}$ & $\sh^{MA_2}$\\ \hline
$c_1$ & 0.2000 & 0.1250 & 0.0822 & 0.1250 & 0.0759\\ \hline
$c_2$ & 0.2000 & 0.1667 & 0.0973 & 0.1200 & 0.0927\\ \hline
$c_3$ & 0.2000 & 0.2500 & 0.1222 & 0.3333 & 0.2310\\ \hline
$c_4$ & 0.2000 & 0.5000 & 0.1746 & 1.5000 & 0.6931\\ \hline
Asy.Risk1 & 3.0752 & 2.0603 &  & 1.5303 & \\ \hline
R.R.R. &  & 33.00 &  & 50.24 & \\ \hline
Asy.Risk2 & 4.0000 &  & 1.8435 &  & 1.4655\\ \hline
R.R.R. &  &  & 53.91 &  & 63.36\\ \hline
\end{tabular}
\begin{tabular}{|l|l|l|l|l|l|} \hline
$n=7$ & $\sh^{U}$ & $\sh^{SDS}$ & $\sh^{KG}$ & $\sh^{MA_1}$ & $\sh^{MA_2}$\\ \hline
$c_1$ & 0.1429 & 0.1000 & 0.0690 & 0.1000 & 0.0647\\ \hline
$c_2$ & 0.1429 & 0.1250 & 0.0796 & 0.0883 & 0.0726\\ \hline
$c_3$ & 0.1429 & 0.1667 & 0.0959 & 0.2000 & 0.1559\\ \hline
$c_4$ & 0.1429 & 0.2500 & 0.1259 & 0.5956 & 0.3816\\ \hline
Asy.Risk1 & 1.8508 & 1.2955 &  & 0.9241 & \\ \hline
R.R.R. &  & 30.01 &  & 50.07 & \\ \hline
Asy.Risk2 & 2.8571 &  & 1.4923 &  & 1.1116\\ \hline
R.R.R. &  &  & 47.77 &  & 61.10\\ \hline
\end{tabular}

\medskip

\begin{tabular}{|l|l|l|l|l|l|} \hline
$n=9$ & $\sh^{U}$ & $\sh^{SDS}$ & $\sh^{KG}$ & $\sh^{MA_1}$ & $\sh^{MA_2}$\\ \hline
$c_1$ & 0.1111 & 0.0833 & 0.0600 & 0.0833 & 0.0571\\ \hline
$c_2$ & 0.1111 & 0.1000 & 0.0681 & 0.0707 & 0.0602\\ \hline
$c_3$ & 0.1111 & 0.1250 & 0.0798 & 0.1429 & 0.1179\\ \hline
$c_4$ & 0.1111 & 0.1667 & 0.0990 & 0.3497 & 0.2553\\ \hline
Asy.Risk1 & 1.3436 & 0.9790 &  & 0.6852 & \\ \hline
R.R.R. &  & 27.13 &  & 49.00 & \\ \hline
Asy.Risk2 & 2.2222 &  & 1.2591 &  & 0.9083\\ \hline
R.R.R. &  &  & 43.34 &  & 59.13\\ \hline
\end{tabular}
\begin{tabular}{|l|l|l|l|l|l|} \hline
$n=11$ & $\sh^{U}$ & $\sh^{SDS}$ & $\sh^{KG}$ & $\sh^{MA_1}$ & $\sh^{MA_2}$\\ \hline
$c_1$ & 0.0909 & 0.0714 & 0.0533 & 0.0714 & 0.0513\\ \hline
$c_2$ & 0.0909 & 0.0833 & 0.0596 & 0.0593 & 0.0517\\ \hline
$c_3$ & 0.0909 & 0.1000 & 0.0685 & 0.1111 & 0.0949\\ \hline
$c_4$ & 0.0909 & 0.1250 & 0.0819 & 0.2413 & 0.1890\\ \hline
Asy.Risk1 & 1.0585 & 0.7956 &  & 0.5496 & \\ \hline
R.R.R. &  & 24.84 &  & 48.08 & \\ \hline
Asy.Risk2 & 1.8182 &  & 1.0927 &  & 0.7730\\ \hline
R.R.R. &  &  & 39.90 &  & 57.49\\ \hline
\end{tabular}

\medskip

\begin{tabular}{|l|l|l|l|l|l|} \hline
$n=21$ & $\sh^{U}$ & $\sh^{SDS}$ & $\sh^{KG}$ & $\sh^{MA_1}$ & $\sh^{MA_2}$\\ \hline
$c_1$ & 0.0476 & 0.0417 & 0.0346 & 0.0417 & 0.0341\\ \hline
$c_2$ & 0.0476 & 0.0455 & 0.0372 & 0.0338 & 0.0311\\ \hline
$c_3$ & 0.0476 & 0.0500 & 0.0404 & 0.0526 & 0.0483\\ \hline
$c_4$ & 0.0476 & 0.0556 & 0.0444 & 0.0879 & 0.0782\\ \hline
Asy.Risk1 & 0.5127 & 0.4183 &  & 0.2769 & \\ \hline
R.R.R. &  & 18.41 &  & 45.99 & \\ \hline
Asy.Risk2 & 0.9524 &  & 0.6708 &  & 0.4526\\ \hline
R.R.R. &  &  & 29.57 &  & 52.47\\ \hline
\end{tabular}
\begin{tabular}{|l|l|l|l|l|l|} \hline
$n=51$ & $\sh^{U}$ & $\sh^{SDS}$ & $\sh^{KG}$ & $\sh^{MA_1}$ & $\sh^{MA_2}$\\ \hline
$c_1$ & 0.0196 & 0.0185 & 0.0170 & 0.0185 & 0.0169\\ \hline
$c_2$ & 0.0196 & 0.0192 & 0.0176 & 0.0154 & 0.0148\\ \hline
$c_3$ & 0.0196 & 0.0200 & 0.0182 & 0.0204 & 0.0197\\ \hline
$c_4$ & 0.0196 & 0.0208 & 0.0189 & 0.0278 & 0.0266\\ \hline
Asy.Risk1 & 0.2016 & 0.1777 &  & 0.1122 & \\ \hline
R.R.R. &  & 11.88 &  & 44.36 & \\ \hline
Asy.Risk2 & 0.3922 &  & 0.3207 &  & 0.2055\\ \hline
R.R.R. &  &  & 18.22 &  & 47.59\\ \hline
\end{tabular}
\end{center}
\end{table} %
\begin{table}[htbp]
\caption{$p=4, m=2$} \label{p=4, m=2}
\begin{center}
\scriptsize
\begin{tabular}{|l|l|l|l|l|l|} \hline
$n=5$ & $\sh^{U}$ & $\sh^{SDS}$ & $\sh^{KG}$ & $\sh^{MA_1}$ & $\sh^{MA_2}$\\ \hline
$c_1$ & 0.2000 & 0.1250 & 0.0822 & 0.1034 & 0.0762\\ \hline
$c_2$ & 0.2000 & 0.1667 & 0.0973 & 0.2308 & 0.1261\\ \hline
$c_3$ & 0.2000 & 0.2500 & 0.1222 & 0.2000 & 0.1173\\ \hline
$c_4$ & 0.2000 & 0.5000 & 0.1746 & 1.0000 & 0.3988\\ \hline
Asy.Risk1 & 3.0752 & 2.2687 &  & 1.9819 & \\ \hline
R.R.R. &  & 26.23 &  & 35.55 & \\ \hline
Asy.Risk2 & 4.0000 &  & 1.8668 &  & 1.7317\\ \hline
R.R.R. &  &  & 53.33 &  & 56.71\\ \hline
\end{tabular}
\begin{tabular}{|l|l|l|l|l|l|} \hline
$n=7$ & $\sh^{U}$ & $\sh^{SDS}$ & $\sh^{KG}$ & $\sh^{MA_1}$ & $\sh^{MA_2}$\\ \hline
$c_1$ & 0.1429 & 0.1000 & 0.0690 & 0.0820 & 0.0632\\ \hline
$c_2$ & 0.1429 & 0.1250 & 0.0796 & 0.1724 & 0.1055\\ \hline
$c_3$ & 0.1429 & 0.1667 & 0.0959 & 0.1304 & 0.0885\\ \hline
$c_4$ & 0.1429 & 0.2500 & 0.1259 & 0.4286 & 0.2425\\ \hline
Asy.Risk1 & 1.8508 & 1.4334 &  & 1.2107 & \\ \hline
R.R.R. &  & 22.55 &  & 34.59 & \\ \hline
Asy.Risk2 & 2.8571 &  & 1.5273 &  & 1.3728\\ \hline
R.R.R. &  &  & 46.54 &  & 51.95\\ \hline
\end{tabular}

\medskip

\begin{tabular}{|l|l|l|l|l|l|} \hline

$n=9$ & $\sh^{U}$ & $\sh^{SDS}$ & $\sh^{KG}$ & $\sh^{MA_1}$ & $\sh^{MA_2}$\\ \hline
$c_1$ & 0.1111 & 0.0833 & 0.0600 & 0.0682 & 0.0546\\ \hline
$c_2$ & 0.1111 & 0.1000 & 0.0681 & 0.1362 & 0.0910\\ \hline
$c_3$ & 0.1111 & 0.1250 & 0.0798 & 0.0980 & 0.0719\\ \hline
$c_4$ & 0.1111 & 0.1667 & 0.0990 & 0.2632 & 0.1727\\ \hline
Asy.Risk1 & 1.3436 & 1.0774 &  & 0.8908 & \\ \hline
R.R.R. &  & 19.81 &  & 33.70 & \\ \hline
Asy.Risk2 & 2.2222 &  & 1.2992 &  & 1.1422\\ \hline
R.R.R. &  &  & 41.54 &  & 48.60\\ \hline
\end{tabular}
\begin{tabular}{|l|l|l|l|l|l|} \hline

$n=11$ & $\sh^{U}$ & $\sh^{SDS}$ & $\sh^{KG}$ & $\sh^{MA_1}$ & $\sh^{MA_2}$\\ \hline
$c_1$ & 0.0909 & 0.0714 & 0.0533 & 0.0586 & 0.0482\\ \hline
$c_2$ & 0.0909 & 0.0833 & 0.0596 & 0.1119 & 0.0798\\ \hline
$c_3$ & 0.0909 & 0.1000 & 0.0685 & 0.0790 & 0.0609\\ \hline
$c_4$ & 0.0909 & 0.1250 & 0.0819 & 0.1872 & 0.1337\\ \hline
Asy.Risk1 & 1.0585 & 0.8700 &  & 0.7080 & \\ \hline
R.R.R. &  & 17.81 &  & 33.11 & \\ \hline
Asy.Risk2 & 1.8182 &  & 1.1337 &  & 0.9792\\ \hline
R.R.R. &  &  & 37.65 &  & 46.14\\ \hline
\end{tabular}

\medskip

\begin{tabular}{|l|l|l|l|l|l|} \hline

$n=21$ & $\sh^{U}$ & $\sh^{SDS}$ & $\sh^{KG}$ & $\sh^{MA_1}$ & $\sh^{MA_2}$\\ \hline
$c_1$ & 0.0476 & 0.0417 & 0.0346 & 0.0349 & 0.0330\\ \hline
$c_2$ & 0.0476 & 0.0455 & 0.0372 & 0.0577 & 0.0531\\ \hline
$c_3$ & 0.0476 & 0.0500 & 0.0404 & 0.0410 & 0.0352\\ \hline
$c_4$ & 0.0476 & 0.0556 & 0.0444 & 0.0735 & 0.0615\\ \hline
Asy.Risk1 & 0.5127 & 0.4477 &  & 0.3477 & \\ \hline
R.R.R. &  & 12.68 &  & 32.18 & \\ \hline
Asy.Risk2 & 0.9524 &  & 0.7013 &  & 0.5722\\ \hline
R.R.R. &  &  & 26.36 &  & 39.92\\ \hline
\end{tabular}
\begin{tabular}{|l|l|l|l|l|l|} \hline

$n=51$ & $\sh^{U}$ & $\sh^{SDS}$ & $\sh^{KG}$ & $\sh^{MA_1}$ & $\sh^{MA_2}$\\ \hline
$c_1$ & 0.0196 & 0.0185 & 0.0170 & 0.0162 & 0.0153\\ \hline
$c_2$ & 0.0196 & 0.0192 & 0.0176 & 0.0227 & 0.0211\\ \hline
$c_3$ & 0.0196 & 0.0200 & 0.0182 & 0.0173 & 0.0163\\ \hline
$c_4$ & 0.0196 & 0.0208 & 0.0189 & 0.0248 & 0.0232\\ \hline
Asy.Risk1 & 0.2016 & 0.1857 &  & 0.1377 & \\ \hline
R.R.R. &  & 7.90 &  & 31.73 & \\ \hline
Asy.Risk2 & 0.3922 &  & 0.3331 &  & 0.2544\\ \hline
R.R.R. &  &  & 15.06 &  & 35.13\\ \hline
\end{tabular}
\end{center}
\end{table} %
\begin{table}[htbp]
\caption{$p=4, m=3$} \label{p=4, m=3}
\begin{center}
\scriptsize
\begin{tabular}{|l|l|l|l|l|l|} \hline

$n=4$ & $\sh^{U}$ & $\sh^{SDS}$ & $\sh^{KG}$ & $\sh^{MA_1}$ & $\sh^{MA_2}$\\ \hline
$c_1$ & 0.2500 & 0.1429 & 0.0919 & 0.1071 & 0.0852\\ \hline
$c_2$ & 0.2500 & 0.2000 & 0.1111 & 0.2500 & 0.1670\\ \hline
$c_3$ & 0.2500 & 0.3333 & 0.1449 & 0.6000 & 0.2383\\ \hline
$c_4$ & 0.2500 & 1.0000 & 0.2174 & 1.0000 & 0.1698\\ \hline
Asy.Risk1 & 4.8592 & 3.6569 &  & 3.4447 & \\ \hline
R.R.R. &  & 24.74 &  & 29.11 & \\ \hline
Asy.Risk2 & 5.0000 &  & 2.0872 &  & 1.9697\\ \hline
R.R.R. &  &  & 58.26 &  & 60.61\\ \hline
\end{tabular}
\begin{tabular}{|l|l|l|l|l|l|} \hline

$n=6$ & $\sh^{U}$ & $\sh^{SDS}$ & $\sh^{KG}$ & $\sh^{MA_1}$ & $\sh^{MA_2}$\\ \hline
$c_1$ & 0.1667 & 0.1111 & 0.0749 & 0.0812 & 0.0678\\ \hline
$c_2$ & 0.1667 & 0.1429 & 0.0873 & 0.1667 & 0.1248\\ \hline
$c_3$ & 0.1667 & 0.2000 & 0.1072 & 0.3733 & 0.2028\\ \hline
$c_4$ & 0.1667 & 0.3333 & 0.1461 & 0.3333 & 0.1209\\ \hline
Asy.Risk1 & 2.2985 & 1.7446 &  & 1.5186 & \\ \hline
R.R.R. &  & 24.10 &  & 33.93 & \\ \hline
Asy.Risk2 & 3.3333 &  & 1.6702 &  & 1.5097\\ \hline
R.R.R. &  &  & 49.89 &  & 54.71\\ \hline
\end{tabular}
\begin{tabular}{|l|l|l|l|l|l|} \hline

$n=8$ & $\sh^{U}$ & $\sh^{SDS}$ & $\sh^{KG}$ & $\sh^{MA_1}$ & $\sh^{MA_2}$\\ \hline
$c_1$ & 0.1250 & 0.0909 & 0.0642 & 0.0660 & 0.0569\\ \hline
$c_2$ & 0.1250 & 0.1111 & 0.0733 & 0.1250 & 0.0999\\ \hline
$c_3$ & 0.1250 & 0.1429 & 0.0870 & 0.2591 & 0.1670\\ \hline
$c_4$ & 0.1250 & 0.2000 & 0.1108 & 0.2000 & 0.0966\\ \hline
Asy.Risk1 & 1.5538 & 1.2032 &  & 0.9929 & \\ \hline
R.R.R. &  & 22.57 &  & 36.10 & \\ \hline
Asy.Risk2 & 2.5000 &  & 1.3948 &  & 1.2111\\ \hline
R.R.R. &  &  & 44.21 &  & 51.56\\ \hline
\end{tabular}
\begin{tabular}{|l|l|l|l|l|l|} \hline

$n=10$ & $\sh^{U}$ & $\sh^{SDS}$ & $\sh^{KG}$ & $\sh^{MA_1}$ & $\sh^{MA_2}$\\ \hline
$c_1$ & 0.1000 & 0.0769 & 0.0565 & 0.0560 & 0.0493\\ \hline
$c_2$ & 0.1000 & 0.0909 & 0.0636 & 0.1000 & 0.0833\\ \hline
$c_3$ & 0.1000 & 0.1111 & 0.0737 & 0.1944 & 0.1385\\ \hline
$c_4$ & 0.1000 & 0.1429 & 0.0896 & 0.1429 & 0.0810\\ \hline
Asy.Risk1 & 1.1828 & 0.9327 &  & 0.7412 & \\ \hline
R.R.R. &  & 21.15 &  & 37.34 & \\ \hline
Asy.Risk2 & 2.0000 &  & 1.1991 &  & 1.0067\\ \hline
R.R.R. &  &  & 40.05 &  & 49.66\\ \hline
\end{tabular}
\begin{tabular}{|l|l|l|l|l|l|} \hline

$n=20$ & $\sh^{U}$ & $\sh^{SDS}$ & $\sh^{KG}$ & $\sh^{MA_1}$ & $\sh^{MA_2}$\\ \hline
$c_1$ & 0.0500 & 0.0435 & 0.0358 & 0.0326 & 0.0303\\ \hline
$c_2$ & 0.0500 & 0.0476 & 0.0386 & 0.0500 & 0.0455\\ \hline
$c_3$ & 0.0500 & 0.0526 & 0.0421 & 0.0808 & 0.0694\\ \hline
$c_4$ & 0.0500 & 0.0588 & 0.0465 & 0.0588 & 0.0450\\ \hline
Asy.Risk1 & 0.5385 & 0.4484 &  & 0.3218 & \\ \hline
R.R.R. &  & 16.73 &  & 40.24 & \\ \hline
Asy.Risk2 & 1.0000 &  & 0.7122 &  & 0.5395\\ \hline
R.R.R. &  &  & 28.78 &  & 46.05\\ \hline
\end{tabular}
\begin{tabular}{|l|l|l|l|l|l|} \hline

$n=50$ & $\sh^{U}$ & $\sh^{SDS}$ & $\sh^{KG}$ & $\sh^{MA_1}$ & $\sh^{MA_2}$\\ \hline
$c_1$ & 0.0200 & 0.0189 & 0.0172 & 0.0151 & 0.0146\\ \hline
$c_2$ & 0.0200 & 0.0196 & 0.0179 & 0.0200 & 0.0192\\ \hline
$c_3$ & 0.0200 & 0.0204 & 0.0186 & 0.0271 & 0.0256\\ \hline
$c_4$ & 0.0200 & 0.0213 & 0.0193 & 0.0213 & 0.0192\\ \hline
Asy.Risk1 & 0.2069 & 0.1836 &  & 0.1207 & \\ \hline
R.R.R. &  & 11.28 &  & 41.65 & \\ \hline
Asy.Risk2 & 0.4000 &  & 0.3293 &  & 0.2236\\ \hline
R.R.R. &  &  & 17.68 &  & 44.11\\ \hline
\end{tabular}
\end{center}
\end{table} %

Tables \ref{p=3, m=1}--\ref{p=4, m=3} give $c_i$'s for $\sh^{U},
\sh^{SDS}, \sh^{KG}, \sh^{MA_1}, \sh^{MA_2}$ when $p=3,4$ with several
values of $n$. The value of $c_i$'s for the minimum asymptotic risk
estimators $\sh^{MA_1}, \sh^{MA_2}$ is calculated by the
aforementioned analytic method. Note that for the case $p=2$, the
minimum asymptotic risk estimator naturally coincides with
$\sh^{SDS}(\sh^{KG})$ which is the minimum asymptotic risk estimator
for $L_1(L_2)$ when we see the total dispersion of population
eigenvalues (see Takemura ans Sheena (2005)). As it is well known,
$n^{-1}l_i \ (i=1,\ldots,p)$ tends to overestimate the corresponding
eigenvalue of $\s$ when $i$ is small, while it tends to underestimate
the corresponding eigenvalue of $\s$ when $i$ is large. The estimators
$\sh^{SDS}$, $\sh^{KG}$ modify this tendency by increasing weight
$c_1<\cdots <c_p$. It is seen from the tables that $\sh^{MA_1},
\sh^{MA_2}$ enlarge the weight difference within each block in most
cases; for example when $p=4, m=2$, the relation between $c_i$'s of
$\sh^{SDS}(\sh^{KG})$ (say $c_i^{SDS}(c_i^{KG}),\ i=1,\ldots,4$) and
those of $\sh^{MA_1}(\sh^{MA_2})$ (say $c_i^{MA_1}(c_i^{MA_2}),\
i=1,\ldots, 4$) is found as
$$
c_1^{MA_1}<c_1^{SDS}<c_2^{SDS}<c_2^{MA_1},\qquad c_3^{MA_1}<c_3^{SDS}<c_4^{SDS}<c_4^{MA_1},
$$
and 
$$
c_1^{MA_2}<c_1^{KG}<c_2^{KG}<c_2^{MA_2},\qquad c_3^{MA_2}<c_3^{KG}<c_4^{KG}<c_4^{MA_2}.
$$

The tables also give asymptotic risk comparison w.r.t.\ $L_1$ among the estimators $\sh^{U}$, $\sh^{SDS}$, $\sh^{MA_1}$ (see ``Asy.Risk1'') and that w.r.t. $L_2$ among the estimators $\sh^{U}, \sh^{KG}, \sh^{MA_2}$ (see ``Asy.Risk2''). The risks are analytically calculated except for evaluating $\sumi E[\log \chi^2_{n-i+1}]$ by Monte Carlo simulation method using $10^5$ random numbers. ``R.R.R.'' under ``Asy.Risk1'' or ``Asy.Risk2'' shows the risk reduction rate defined by
$$
\mbox{R.R.R. of }\sh=\frac{\mbox{The risk of }\sh^{U}-\mbox{The risk of }\sh}{\mbox{The risk of }\sh^{U}}
\times 100.
$$
It has been observed that $\sh^{SDS}$ and $\sh^{KG}$ drastically reduce the risk of $\s^{U}$ when the population eigenvalues are close to each other. Lin and Perlman (1985) reports that when $\s=\ii_p$,  R.R.R. of $\sh^{SDS}$ often reaches 70\%. See also Sugiura and Ishibayashi (1997) for a risk comparison by elabarate simulation. In the situation of the block-wise dispersion, the risk reduction rate of these estimators rarely approaches 50\%. Especially when $n$ is as large as 50, the rate is always under 20\%. On the other hand, the risk reduction rates of $\sh^{MA_1}$ and $\s^{MA_2}$ are constantly over 30\% and often reach 50\% irrespective of the values of $n$. It is interesting that $\sh^{MA_2}$ always outperforms $\sh^{MA_1}$ in view of R.R.R.
\subsection{Simulation studies}
\label{subsec:numerical_analysis}
In this subsection, we evaluate the performance of $\sh^{MA_1},\
\sh^{MA_2}$ by Monte Carlo simulation
under the situation (\ref{pract_case}). 
As we saw in the previous
subsection, in view of the asymptotic risks, $\sh^{MA_1},\ \sh^{MA_2}$
provide better risk reduction compared to $\sh^{SDS},\ \sh^{KG}$. In
practical point view, however, it is important to see how largely the
population eigenvalues must be dispersed so that the use of
$\sh^{MA_\lossd},\ \lossd=1,2$, is recommended. The convergence speed of the
distributions given in Theorem \ref{thm:asympt_dist}, which is an interesting topic by itself, is closely related to this problem.

To see the convergence speed in both distributions and risks, we
carried out Monte Carlo Simulation for the two cases $p=3,\ m=1$ and
$p=4,\ m=1$. In each case, we took 11 values 1.0, 0.8, 0.6, 0.4, 0.2,
$10^{-i}(i=1,\ldots,6)$ in the convergence parameter $\beta$, while $\alpha$ is
fixed at 1.
We took three different values
of $n$ in each case and generated $10^6$ random Wishart matrices under
given $p, n, \beta$. The result is given in Table \ref{speed p=3, m=1}
($p=3, m=1$) and Table \ref{speed p=4, m=1} ($p=4, m=1$). The upper part of
each table shows the speed of the distributional convergence in
Theorem \ref{thm:asympt_dist}. Note that when $\xxxi_1=\ii_m,\
\xxxi_2=\ii_{p-m},$ the asymptotic distribution of a diagonal element
of $\widetilde{\ww}_{ss},\ s=1,2$, is a $\chi^2$ distribution. The labels in the
tables are given as follows with $\chi^2_n(\alpha),\ z(\alpha)$
denoting the lower $\alpha$ percentage points of $\chi^2$ distribution
with $n$ degrees of freedom and the standard normal distribution, respectively ;

$$
\begin{array}{rll}
&\mbox{Table 6}&\\
&\mbox{Prob 1a}=P(\widetilde{\ww}_{11}\leq \chi^2_n(0.05)),\quad &\mbox{Prob 1b}=P(\widetilde{\ww}_{11}\leq \chi^2_n(0.95)),\\
&\mbox{Prob 2a}=P((\widetilde{\ww}_{22})_{11}\leq \chi^2_{n-1}(0.05)),\quad &\mbox{Prob 2b}=P((\widetilde{\ww}_{22})_{11}\leq \chi^2_{n-1}(0.95)),\\
&\mbox{Prob 3a}=P((\widetilde{\ww}_{22})_{22}\leq \chi^2_{n-1}(0.05)),\quad &\mbox{Prob 3b}=P((\widetilde{\ww}_{22})_{22}\leq \chi^2_{n-1}(0.95)),\\
&\mbox{Prob 4a}=P((\widetilde{{\bm Z}}_{21})_{11}\leq z(0.05)),\quad &\mbox{Prob 4b}=P((\widetilde{{\bm Z}}_{21})_{11}\leq z(0.95)),\\
&\mbox{Prob 5a}=P((\widetilde{{\bm Z}}_{21})_{21}\leq z(0.05)),\quad &\mbox{Prob 5b}=P((\widetilde{{\bm Z}}_{21})_{21}\leq z(0.95)),\\
\end{array}
$$
$$
\begin{array}{rll}
&\mbox{Table 7}&\\
&\mbox{Prob 1a}=P(\widetilde{\ww}_{11}\leq \chi^2_n(0.05)),\quad &\mbox{Prob 1b}=P(\widetilde{\ww}_{11}\leq \chi^2_n(0.95)),\\
&\mbox{Prob 2a}=P((\widetilde{\ww}_{22})_{11}\leq \chi^2_{n-1}(0.05)),\quad &\mbox{Prob 2b}=P((\widetilde{\ww}_{22})_{11}\leq \chi^2_{n-1}(0.95)),\\
&\mbox{Prob 3a}=P((\widetilde{\ww}_{22})_{33}\leq \chi^2_{n-1}(0.05)),\quad &\mbox{Prob 3b}=P((\widetilde{\ww}_{22})_{33}\leq \chi^2_{n-1}(0.95)),\\
&\mbox{Prob 4a}=P((\widetilde{{\bm Z}}_{21})_{11}\leq z(0.05)),\quad &\mbox{Prob 4b}=P((\widetilde{{\bm Z}}_{21})_{11}\leq z(0.95)),\\
&\mbox{Prob 5a}=P((\widetilde{{\bm Z}}_{21})_{31}\leq z(0.05)),\quad &\mbox{Prob 5b}=P((\widetilde{{\bm Z}}_{21})_{31}\leq z(0.95)).
\end{array}
$$
\begin{table}[tbp]
\caption{$p=3, m=1$} \label{speed p=3, m=1}
\begin{center}
\scriptsize
\begin{tabular}{|l|l|l|l|l|l|l|l|l|l|l|l|l|} \hline
$n=10$ & 1 & 0.8 & 0.6 & 0.4 & 0.2 & $10^{-1}$ & $10^{-2}$ & $10^{-3}$ & $10^{-4}$ & $10^{-5}$ & $10^{-6}$ & Asymp.\\ \hline
Prob 1a & 0.4994 & 0.3992 & 0.2814 & 0.1551 & 0.0695 & 0.0534 & 0.0501 & 0.0491 & 0.0507 & 0.0491 & 0.0508 & 0.0500\\ \hline
Prob 2a & 0.4091 & 0.3273 & 0.2321 & 0.1321 & 0.0677 & 0.0558 & 0.0516 & 0.0489 & 0.0504 & 0.0495 & 0.0502 & 0.0500\\ \hline
Prob 3a & 0.4121 & 0.3302 & 0.2311 & 0.1317 & 0.0684 & 0.0564 & 0.0503 & 0.0499 & 0.0505 & 0.0499 & 0.0518 & 0.0500\\ \hline
Prob 4a & 0.2024 & 0.1799 & 0.1502 & 0.1072 & 0.0597 & 0.0385 & 0.0263 & 0.0255 & 0.0294 & 0.0429 & 0.0499 & 0.0500\\ \hline
Prob 5a & 0.0001 & 0.0000 & 0.0000 & 0.0000 & 0.0000 & 0.0000 & 0.0000 & 0.0000 & 0.0000 & 0.0269 & 0.0500 & 0.0500\\ \hline
Prob 1b & 0.9700 & 0.9636 & 0.9572 & 0.9531 & 0.9503 & 0.9507 & 0.9499 & 0.9514 & 0.9496 & 0.9497 & 0.9496 & 0.9500\\ \hline
Prob 2b & 0.9993 & 0.9985 & 0.9955 & 0.9871 & 0.9695 & 0.9576 & 0.9508 & 0.9498 & 0.9503 & 0.9488 & 0.9498 & 0.9500\\ \hline
Prob 3b & 0.9994 & 0.9983 & 0.9957 & 0.9874 & 0.9693 & 0.9582 & 0.9508 & 0.9509 & 0.9496 & 0.9500 & 0.9497 & 0.9500\\ \hline
Prob 4b & 0.6174 & 0.6492 & 0.6986 & 0.7671 & 0.8528 & 0.8924 & 0.9236 & 0.9255 & 0.9301 & 0.9451 & 0.9515 & 0.9500\\ \hline
Prob 5b & 0.4137 & 0.4673 & 0.5465 & 0.6624 & 0.7937 & 0.8530 & 0.8971 & 0.8994 & 0.9001 & 0.9275 & 0.9504 & 0.9500\\ \hline
Risk 1\_U & 0.6769 & 0.6753 & 0.6786 & 0.6779 & 0.6777 & 0.6778 & 0.6784 & 0.6757 & 0.6759 & 0.6758 & 0.6800 & 0.6765\\ \hline
Risk 1\_SDS & 0.4589 & 0.4611 & 0.4770 & 0.5038 & 0.5409 & 0.5580 & 0.5701 & 0.5687 & 0.5690 & 0.5684 & 0.5727 & 0.5692\\ \hline
Risk 1\_MA1 & 0.3595 & 0.3644 & 0.3824 & 0.4091 & 0.4400 & 0.4553 & 0.4677 & 0.4668 & 0.4677 & 0.4660 & 0.4704 & 0.4676\\ \hline
Risk 2\_U & 1.1996 & 1.1997 & 1.2017 & 1.1976 & 1.1983 & 1.1989 & 1.1980 & 1.1966 & 1.2020 & 1.1990 & 1.2021 & 1.2000\\ \hline
Risk 2\_KG & 0.7117 & 0.7132 & 0.7228 & 0.7407 & 0.7641 & 0.7748 & 0.7815 & 0.7812 & 0.7806 & 0.7811 & 0.7839 & 0.7817\\ \hline
Risk 2\_MA2 & 0.6109 & 0.6147 & 0.6255 & 0.6397 & 0.6540 & 0.6625 & 0.6706 & 0.6703 & 0.6704 & 0.6689 & 0.6725 & 0.6591\\ \hline
\end{tabular}

\medskip

\begin{tabular}{|l|l|l|l|l|l|l|l|l|l|l|l|l|} \hline
$n=20$ & 1 & 0.8 & 0.6 & 0.4 & 0.2 & $10^{-1}$ & $10^{-2}$ & $10^{-3}$ & $10^{-4}$ & $10^{-5}$ & $10^{-6}$ & Asymp.\\ \hline
Prob 1a & 0.6164 & 0.4550 & 0.2706 & 0.1187 & 0.0574 & 0.0523 & 0.0495 & 0.0506 & 0.0517 & 0.0491 & 0.0498 & 0.0500\\ \hline
Prob 2a & 0.5114 & 0.3837 & 0.2285 & 0.1027 & 0.0608 & 0.0537 & 0.0511 & 0.0498 & 0.0505 & 0.0490 & 0.0494 & 0.0500\\ \hline
Prob 3a & 0.5081 & 0.3856 & 0.2309 & 0.1043 & 0.0594 & 0.0528 & 0.0511 & 0.0508 & 0.0493 & 0.0513 & 0.0499 & 0.0500\\ \hline
Prob 4a & 0.2493 & 0.2196 & 0.1684 & 0.1100 & 0.0560 & 0.0377 & 0.0257 & 0.0264 & 0.0328 & 0.0483 & 0.0513 & 0.0500\\ \hline
Prob 5a & 0.0015 & 0.0015 & 0.0008 & 0.0002 & 0.0000 & 0.0000 & 0.0000 & 0.0000 & 0.0003 & 0.0451 & 0.0497 & 0.0500\\ \hline
Prob 1b & 0.9767 & 0.9661 & 0.9595 & 0.9543 & 0.9515 & 0.9484 & 0.9513 & 0.9500 & 0.9498 & 0.9498 & 0.9494 & 0.9500\\ \hline
Prob 2b & 0.9995 & 0.9984 & 0.9936 & 0.9816 & 0.9631 & 0.9547 & 0.9513 & 0.9499 & 0.9498 & 0.9511 & 0.9500 & 0.9500\\ \hline
Prob 3b & 0.9994 & 0.9983 & 0.9940 & 0.9815 & 0.9623 & 0.9552 & 0.9520 & 0.9516 & 0.9499 & 0.9501 & 0.9505 & 0.9500\\ \hline
Prob 4b & 0.5542 & 0.6008 & 0.6689 & 0.7653 & 0.8592 & 0.8975 & 0.9200 & 0.9257 & 0.9334 & 0.9499 & 0.9508 & 0.9500\\ \hline
Prob 5b & 0.3123 & 0.3793 & 0.4993 & 0.6566 & 0.8026 & 0.8574 & 0.8953 & 0.8998 & 0.8987 & 0.9457 & 0.9503 & 0.9500\\ \hline
Risk 1\_U & 0.3178 & 0.3179 & 0.3181 & 0.3171 & 0.3183 & 0.3178 & 0.3177 & 0.3175 & 0.3177 & 0.3171 & 0.3176 & 0.3164\\ \hline
Risk 1\_SDS & 0.2363 & 0.2390 & 0.2486 & 0.2628 & 0.2767 & 0.2802 & 0.2829 & 0.2830 & 0.2833 & 0.2827 & 0.2832 & 0.2819\\ \hline
Risk 1\_MA1 & 0.1880 & 0.1923 & 0.2023 & 0.2115 & 0.2200 & 0.2236 & 0.2261 & 0.2262 & 0.2267 & 0.2260 & 0.2262 & 0.2251\\ \hline
Risk 2\_U & 0.5995 & 0.6011 & 0.6008 & 0.5992 & 0.5999 & 0.6006 & 0.5992 & 0.5987 & 0.6005 & 0.6003 & 0.6011 & 0.6000\\ \hline
Risk 2\_KG & 0.4085 & 0.4117 & 0.4226 & 0.4384 & 0.4530 & 0.4568 & 0.4595 & 0.4598 & 0.4600 & 0.4593 & 0.4594 & 0.4598\\ \hline
Risk 2\_MA2 & 0.3563 & 0.3606 & 0.3674 & 0.3706 & 0.3744 & 0.3775 & 0.3792 & 0.3794 & 0.3801 & 0.3797 & 0.3793 & 0.3745\\ \hline
\end{tabular}

\medskip

\begin{tabular}{|l|l|l|l|l|l|l|l|l|l|l|l|l|} \hline
$n=50$ & 1 & 0.8 & 0.6 & 0.4 & 0.2 & $10^{-1}$ & $10^{-2}$ & $10^{-3}$ & $10^{-4}$ & $10^{-5}$ & $10^{-6}$ & Asymp.\\ \hline
Prob 1a & 0.7358 & 0.4769 & 0.2076 & 0.0793 & 0.0532 & 0.0503 & 0.0484 & 0.0506 & 0.0506 & 0.0485 & 0.0501 & 0.0500\\ \hline
Prob 2a & 0.6110 & 0.4089 & 0.1725 & 0.0720 & 0.0549 & 0.0511 & 0.0512 & 0.0498 & 0.0505 & 0.0489 & 0.0499 & 0.0500\\ \hline
Prob 3a & 0.6079 & 0.4100 & 0.1737 & 0.0732 & 0.0564 & 0.0513 & 0.0495 & 0.0484 & 0.0487 & 0.0489 & 0.0504 & 0.0500\\ \hline
Prob 4a & 0.2992 & 0.2529 & 0.1788 & 0.1042 & 0.0541 & 0.0368 & 0.0271 & 0.0274 & 0.0411 & 0.0493 & 0.0500 & 0.0500\\ \hline
Prob 5a & 0.0200 & 0.0109 & 0.0031 & 0.0002 & 0.0000 & 0.0000 & 0.0000 & 0.0000 & 0.0072 & 0.0490 & 0.0506 & 0.0500\\ \hline
Prob 1b & 0.9823 & 0.9707 & 0.9606 & 0.9532 & 0.9511 & 0.9499 & 0.9485 & 0.9503 & 0.9489 & 0.9500 & 0.9505 & 0.9500\\ \hline
Prob 2b & 0.9995 & 0.9979 & 0.9883 & 0.9696 & 0.9567 & 0.9519 & 0.9503 & 0.9511 & 0.9498 & 0.9506 & 0.9498 & 0.9500\\ \hline
Prob 3b & 0.9996 & 0.9977 & 0.9889 & 0.9705 & 0.9557 & 0.9528 & 0.9495 & 0.9499 & 0.9500 & 0.9506 & 0.9498 & 0.9500\\ \hline
Prob 4b & 0.5050 & 0.5599 & 0.6582 & 0.7701 & 0.8610 & 0.8967 & 0.9231 & 0.9266 & 0.9408 & 0.9502 & 0.9495 & 0.9500\\ \hline
Prob 5b & 0.2273 & 0.3156 & 0.4805 & 0.6648 & 0.8100 & 0.8613 & 0.8960 & 0.8993 & 0.9066 & 0.9494 & 0.9498 & 0.9500\\ \hline
Risk 1\_U & 0.1223 & 0.1226 & 0.1227 & 0.1229 & 0.1228 & 0.1228 & 0.1226 & 0.1223 & 0.1228 & 0.1221 & 0.1230 & 0.1236\\ \hline
Risk 1\_SDS & 0.1006 & 0.1026 & 0.1074 & 0.1117 & 0.1137 & 0.1143 & 0.1145 & 0.1143 & 0.1148 & 0.1140 & 0.1149 & 0.1155\\ \hline
Risk 1MA1 & 0.0814 & 0.0843 & 0.0867 & 0.0871 & 0.0882 & 0.0888 & 0.0891 & 0.0891 & 0.0896 & 0.0887 & 0.0895 & 0.0901\\ \hline
Risk 2\_U & 0.2391 & 0.2399 & 0.2400 & 0.2401 & 0.2402 & 0.2404 & 0.2403 & 0.2391 & 0.2405 & 0.2389 & 0.2406 & 0.2400\\ \hline
Risk 2\_KG & 0.1874 & 0.1906 & 0.1982 & 0.2049 & 0.2079 & 0.2087 & 0.2091 & 0.2086 & 0.2096 & 0.2083 & 0.2097 & 0.2093\\ \hline
Risk 2\_MA2 & 0.1641 & 0.1673 & 0.1669 & 0.1643 & 0.1649 & 0.1654 & 0.1658 & 0.1656 & 0.1665 & 0.1650 & 0.1663 & 0.1647\\ \hline
\end{tabular}
\end{center}
\end{table} %
\begin{table}[tbp]
\caption{$p=4, m=1$} \label{speed p=4, m=1}
\begin{center}
\scriptsize
\begin{tabular}{|l|l|l|l|l|l|l|l|l|l|l|l|l|} \hline
$n=11$ & 1 & 0.8 & 0.6 & 0.4 & 0.2 & $10^{-1}$ & $10^{-2}$ & $10^{-3}$ & $10^{-4}$ & $10^{-5}$ & $10^{-6}$ & Asymp.\\ \hline
Prob 1a & 0.5852 & 0.4760 & 0.3396 & 0.1881 & 0.0751 & 0.0554 & 0.0517 & 0.0497 & 0.0495 & 0.0498 & 0.0489 & 0.0500\\ \hline
Prob 2a & 0.3156 & 0.2620 & 0.1937 & 0.1166 & 0.0660 & 0.0557 & 0.0505 & 0.0494 & 0.0510 & 0.0501 & 0.0492 & 0.0500\\ \hline
Prob 3a & 0.3146 & 0.2620 & 0.1920 & 0.1172 & 0.0651 & 0.0560 & 0.0515 & 0.0503 & 0.0496 & 0.0494 & 0.0503 & 0.0500\\ \hline
Prob 4a & 0.1955 & 0.1794 & 0.1533 & 0.1124 & 0.0650 & 0.0431 & 0.0283 & 0.0284 & 0.0315 & 0.0449 & 0.0507 & 0.0500\\ \hline
Prob 5a & 0.1312 & 0.1206 & 0.1005 & 0.0685 & 0.0342 & 0.0205 & 0.0131 & 0.0131 & 0.0172 & 0.0391 & 0.0497 & 0.0500\\ \hline
Prob 1b & 0.9761 & 0.9676 & 0.9610 & 0.9547 & 0.9521 & 0.9510 & 0.9500 & 0.9510 & 0.9495 & 0.9508 & 0.9501 & 0.9500\\ \hline
Prob 2b & 0.9986 & 0.9977 & 0.9948 & 0.9863 & 0.9690 & 0.9581 & 0.9505 & 0.9510 & 0.9507 & 0.9509 & 0.9493 & 0.9500\\ \hline
Prob 3b & 0.9985 & 0.9979 & 0.9947 & 0.9861 & 0.9681 & 0.9570 & 0.9509 & 0.9515 & 0.9499 & 0.9508 & 0.9508 & 0.9500\\ \hline
Prob 4b & 0.6525 & 0.6774 & 0.7152 & 0.7789 & 0.8576 & 0.8979 & 0.9251 & 0.9276 & 0.9310 & 0.9454 & 0.9501 & 0.9500\\ \hline
Prob 5b & 0.5899 & 0.6184 & 0.6607 & 0.7327 & 0.8304 & 0.8751 & 0.9091 & 0.9115 & 0.9185 & 0.9399 & 0.9508 & 0.9500\\ \hline
Risk 1\_U & 1.0566 & 1.0583 & 1.0552 & 1.0592 & 1.0577 & 1.0583 & 1.0603 & 1.0544 & 1.0574 & 1.0573 & 1.0559 & 1.0585\\ \hline
Risk 1\_SDS & 0.6514 & 0.6572 & 0.6714 & 0.7092 & 0.7558 & 0.7781 & 0.7954 & 0.7920 & 0.7943 & 0.7942 & 0.7927 & 0.7956\\ \hline
Risk 1\_MA1 & 0.4064 & 0.4154 & 0.4367 & 0.4738 & 0.5104 & 0.5295 & 0.5485 & 0.5471 & 0.5484 & 0.5478 & 0.5468 & 0.5496\\ \hline
Risk 2\_U & 1.8199 & 1.8213 & 1.8147 & 1.8170 & 1.8175 & 1.8199 & 1.8210 & 1.8147 & 1.8206 & 1.8180 & 1.8176 & 1.8182\\ \hline
Risk 2\_KG & 1.0173 & 1.0214 & 1.0291 & 1.0493 & 1.0749 & 1.0876 & 1.0939 & 1.0921 & 1.0929 & 1.0926 & 1.0915 & 1.0927\\ \hline
Risk 2\_MA2 & 0.5967 & 0.6075 & 0.6326 & 0.6767 & 0.7268 & 0.7516 & 0.7728 & 0.7724 & 0.7737 & 0.7733 & 0.7719 & 0.7730\\ \hline
\end{tabular}

\medskip

\begin{tabular}{|l|l|l|l|l|l|l|l|l|l|l|l|l|} \hline
$n=21$ & 1 & 0.8 & 0.6 & 0.4 & 0.2 & $10^{-1}$ & $10^{-2}$ & $10^{-3}$ & $10^{-4}$ & $10^{-5}$ & $10^{-6}$ & Asymp.\\ \hline
Prob 1a & 0.7030 & 0.5419 & 0.3317 & 0.1428 & 0.0601 & 0.0532 & 0.0503 & 0.0505 & 0.0498 & 0.0495 & 0.0509 & 0.0500\\ \hline
Prob 2a & 0.4043 & 0.3183 & 0.2017 & 0.0985 & 0.0579 & 0.0547 & 0.0508 & 0.0495 & 0.0492 & 0.0505 & 0.0505 & 0.0500\\ \hline
Prob 3a & 0.3995 & 0.3222 & 0.2019 & 0.0975 & 0.0581 & 0.0527 & 0.0507 & 0.0504 & 0.0496 & 0.0497 & 0.0493 & 0.0500\\ \hline
Prob 4a & 0.2413 & 0.2156 & 0.1748 & 0.1172 & 0.0601 & 0.0421 & 0.0297 & 0.0292 & 0.0344 & 0.0480 & 0.0503 & 0.0500\\ \hline
Prob 5a & 0.1720 & 0.1533 & 0.1185 & 0.0711 & 0.0331 & 0.0201 & 0.0141 & 0.0137 & 0.0211 & 0.0484 & 0.0505 & 0.0500\\ \hline
Prob 1b & 0.9830 & 0.9737 & 0.9627 & 0.9557 & 0.9502 & 0.9506 & 0.9503 & 0.9504 & 0.9497 & 0.9505 & 0.9505 & 0.9500\\ \hline
Prob 2b & 0.9989 & 0.9977 & 0.9929 & 0.9809 & 0.9617 & 0.9548 & 0.9507 & 0.9501 & 0.9488 & 0.9500 & 0.9514 & 0.9500\\ \hline
Prob 3b & 0.9988 & 0.9975 & 0.9935 & 0.9805 & 0.9628 & 0.9549 & 0.9509 & 0.9502 & 0.9501 & 0.9496 & 0.9487 & 0.9500\\ \hline
Prob 4b & 0.5985 & 0.6278 & 0.6881 & 0.7757 & 0.8657 & 0.8998 & 0.9262 & 0.9272 & 0.9339 & 0.9476 & 0.9494 & 0.9500\\ \hline
Prob 5b & 0.5303 & 0.5632 & 0.6291 & 0.7282 & 0.8368 & 0.8794 & 0.9118 & 0.9130 & 0.9219 & 0.9479 & 0.9504 & 0.9500\\ \hline
Risk 1\_U & 0.5121 & 0.5136 & 0.5135 & 0.5116 & 0.5115 & 0.5128 & 0.5110 & 0.5127 & 0.5115 & 0.5109 & 0.5119 & 0.5127\\ \hline
Risk 1\_SDS & 0.3503 & 0.3552 & 0.3677 & 0.3871 & 0.4056 & 0.4134 & 0.4167 & 0.4183 & 0.4172 & 0.4169 & 0.4177 & 0.4183\\ \hline
Risk 1\_MA1 & 0.2241 & 0.2315 & 0.2461 & 0.2568 & 0.2650 & 0.2715 & 0.2759 & 0.2772 & 0.2759 & 0.2764 & 0.2765 & 0.2769\\ \hline
Risk 2\_U & 0.9512 & 0.9514 & 0.9537 & 0.9503 & 0.9516 & 0.9521 & 0.9477 & 0.9535 & 0.9520 & 0.9501 & 0.9505 & 0.9524\\ \hline
Risk 2\_KG & 0.6059 & 0.6109 & 0.6233 & 0.6429 & 0.6607 & 0.6669 & 0.6692 & 0.6708 & 0.6700 & 0.6695 & 0.6707 & 0.6708\\ \hline
Risk 2\_MA2 & 0.3510 & 0.3622 & 0.3861 & 0.4114 & 0.4326 & 0.4433 & 0.4516 & 0.4532 & 0.4521 & 0.4524 & 0.4525 & 0.4526\\ \hline
\end{tabular}

\medskip

\begin{tabular}{|l|l|l|l|l|l|l|l|l|l|l|l|l|} \hline
$n=51$ & 1 & 0.8 & 0.6 & 0.4 & 0.2 & $10^{-1}$ & $10^{-2}$ & $10^{-3}$ & $10^{-4}$ & $10^{-5}$ & $10^{-6}$ & Asymp.\\ \hline
Prob 1a & 0.8209 & 0.5805 & 0.2691 & 0.0916 & 0.0533 & 0.0492 & 0.0498 & 0.0504 & 0.0501 & 0.0504 & 0.0500 & 0.0500\\ \hline
Prob 2a & 0.5101 & 0.3626 & 0.1647 & 0.0721 & 0.0560 & 0.0522 & 0.0501 & 0.0502 & 0.0501 & 0.0504 & 0.0500 & 0.0500\\ \hline
Prob 3a & 0.5098 & 0.3610 & 0.1669 & 0.0722 & 0.0555 & 0.0533 & 0.0500 & 0.0507 & 0.0479 & 0.0498 & 0.0506 & 0.0500\\ \hline
Prob 4a & 0.2912 & 0.2595 & 0.1878 & 0.1118 & 0.0604 & 0.0415 & 0.0303 & 0.0291 & 0.0403 & 0.0507 & 0.0491 & 0.0500\\ \hline
Prob 5a & 0.2191 & 0.1863 & 0.1307 & 0.0689 & 0.0308 & 0.0196 & 0.0133 & 0.0148 & 0.0313 & 0.0501 & 0.0499 & 0.0500\\ \hline
Prob 1b & 0.9891 & 0.9762 & 0.9649 & 0.9548 & 0.9507 & 0.9501 & 0.9504 & 0.9501 & 0.9504 & 0.9505 & 0.9497 & 0.9500\\ \hline
Prob 2b & 0.9992 & 0.9970 & 0.9889 & 0.9700 & 0.9573 & 0.9526 & 0.9502 & 0.9498 & 0.9503 & 0.9507 & 0.9504 & 0.9500\\ \hline
Prob 3b & 0.9990 & 0.9973 & 0.9891 & 0.9712 & 0.9565 & 0.9521 & 0.9499 & 0.9513 & 0.9503 & 0.9506 & 0.9492 & 0.9500\\ \hline
Prob 4b & 0.5383 & 0.5836 & 0.6703 & 0.7803 & 0.8683 & 0.9022 & 0.9272 & 0.9286 & 0.9411 & 0.9494 & 0.9503 & 0.9500\\ \hline
Prob 5b & 0.4666 & 0.5129 & 0.6101 & 0.7334 & 0.8386 & 0.8789 & 0.9081 & 0.9153 & 0.9312 & 0.9496 & 0.9503 & 0.9500\\ \hline
Risk 1\_U & 0.2018 & 0.2022 & 0.2019 & 0.2017 & 0.2019 & 0.2017 & 0.2020 & 0.2017 & 0.2020 & 0.2023 & 0.2018 & 0.2016\\ \hline
Risk 1\_SDS & 0.1566 & 0.1592 & 0.1658 & 0.1721 & 0.1758 & 0.1768 & 0.1780 & 0.1777 & 0.1780 & 0.1783 & 0.1779 & 0.1777\\ \hline
Risk 1MA1 & 0.1037 & 0.1083 & 0.1109 & 0.1088 & 0.1104 & 0.1113 & 0.1125 & 0.1124 & 0.1123 & 0.1124 & 0.1125 & 0.1122\\ \hline
Risk 2\_U & 0.3923 & 0.3939 & 0.3920 & 0.3916 & 0.3920 & 0.3924 & 0.3927 & 0.3919 & 0.3931 & 0.3929 & 0.3920 & 0.3922\\ \hline
Risk 2\_KG & 0.2896 & 0.2938 & 0.3038 & 0.3127 & 0.3179 & 0.3194 & 0.3208 & 0.3203 & 0.3211 & 0.3215 & 0.3208 & 0.3207\\ \hline
Risk 2\_MA2 & 0.1785 & 0.1867 & 0.1943 & 0.1959 & 0.2010 & 0.2033 & 0.2057 & 0.2055 & 0.2054 & 0.2056 & 0.2059 & 0.2055\\ \hline
\end{tabular}
\end{center}
\end{table} %
In the lower part of each table, ``Risk 1\_*'' and ``Risk 2\_*'' show the risks of the corresponding estimator $\sh^*$ respectively for $L_1$ and $L_2$.
The tables show that
\begin{enumerate}
\item The convergence of the diagonal elements of $\widetilde{\ww}_{ss},\ s=1,2$,
  is so rapid that when $\beta=0.1$, the asymptotic distribution already
  gives a good approximation for the exact distribution. When $\beta=0.1$,
  every probability of the diagonal elements is within 0.01 deviation from the exact asymptotic probability.
\item The convergence speed of $\widetilde{\bm{Z}}$ is quite slow compared to that of the diagonal elements of $\widetilde{\ww}_{ss},\ s=1,2$. For a good approximation as above, $\beta$ must be as small as $10^{-5}$ or $10^{-6}$.
\item The risks also rapidly converge to the asymptotic risks so that $\beta=0.1$ is small enough to give a good approximation. Actually all the risks in the tables when $\beta=0.1$ are within the $\pm5\%$ interval centered at the exact asymptotic risk.
\item The risk of $\sh^{MA_\lossd}, \lossd=1,2$, is always lower than
  that of the competing estimators. Most notably their superiority in
  risk is kept even when the population eigenvalues are all equal. It
  seems that $\sh^{MA_\lossd},\ \lossd=1,2$, has robustness to the deviation
  from the dispersion of the population eigenvalues.
\end{enumerate}

Because of the robustness, $\sh^{MA_\lossd},\ \lossd=1,2$, seem to be useful for
various applications.  Now as the last topic in this section, apart
from a decision-theoretic approach, we evaluate these new estimators'
performance in discriminant analysis. We use a well-known example
of Fisher's iris data. The data consists of 50 samples from each of
the three groups(species) with 4-dimensional variable ($x_1$:sepal
length(cm), $x_2$:sepal width(cm), $x_3$:petal length(cm), $x_4$:petal
width(cm)). We downloaded the data from the website {\tt
  http://www-unix.oit.umass.edu/\~{}statdata}.
We let
$\bm{x}_j^{(i)},\ i=1,2,3,\ j=1,\ldots,50$ denote the $j$th
sample in the $i$th group. The estimator to be tested are the
traditional estimators $\sh^{U}$, $\sh^{SDS}$, $\sh^{KG}$ and the new
estimators $\sh^{MA_1}$, $\sh^{MA_2}$ which are formulated under the
condition $p=4, m=1.$

We carry out cross validations. Suppose a learning data set
$\bm{y}^{(i)}_j,\ j=1,\ldots, N$, is chosen from the $i$th group,
$i=1,2,3$. Estimates for the population covariance matrix of the $i$th group 
are calculated from $\sh^{U}$, $\sh^{SDS}$, $\sh^{KG}$, $\sh^{MA_1}$,
$\sh^{MA_2}$ %
based on
$$
\bm{A}^{(i)}=\sum_{j=1}^N (\bm{y}^{(i)}_j-\bar{\bm{y}}^{(i)})(\bm{y}^{(i)}_j-\bar{\bm{y}}^{(i)})',
$$
where $\bar{\bm{y}}^{(i)}=N^{-1}\sum_{j=1}^N\bm{y}_j^{(i)}$.
As a discriminant function, we use a Mahalanobis distance based on each estimates $\sh^U(\bm{A}^{(i)})$, $\sh^{SDS}(\bm{A}^{(i)})$, $\sh^{KG}(\bm{A}^{(i)})$, $\sh^{MA_1}(\bm{A}^{(i)})$, $\sh^{MA_2}(\bm{A}^{(i)})$, that is, for a test data $\bm{x}$ 
$$
MD_i^{*}=(\bm{x}-\bar{\bm{y}}^{(i)})'\sh^{*}(\bm{A}^{(i)})^{-1}(\bm{x}-\bar{\bm{y}}^{(i)}),\quad i=1,2,3.
$$

The eigenvalues of the covariance matrix within each group is as follows; 
\begin{equation}
\label{iris_data_eigenvalues}
\begin{array}{rl}
\mbox{Group 1:}&(0.234, 0.039, 0.027, 0.009),\\ 
\mbox{Group 2:}&(0.482, 0.075, 0.056, 0.011),\\ 
\mbox{Group 3:}&(0.688, 0.107, 0.057, 0.036).
\end{array}
\end{equation}
\begin{table}[tbp]
\caption{10-sample-set} \label{10-sample-set}
\footnotesize
\begin{center}
\begin{tabular}{|c|c|c|c|c|c|} \hline
\multicolumn{1}{|p{40pt}|}{Learning Data Set} & \multicolumn{1}{p{25pt}|}{$\sh^{U}$} & \multicolumn{1}{p{25pt}|}{$\sh^{SDS}$} & \multicolumn{1}{p{25pt}|}{$\sh^{KG}$} & \multicolumn{1}{p{25pt}|}{$\sh^{MA_1}$} & \multicolumn{1}{p{25pt}|}{$\sh^{MA_2}$}\\ \hline
1 & 82.50 & 83.33 & 83.33 & 81.67 & 82.50 \\ \hline
2 & 85.83 & 85.00 & 85.00 & 85.00 & 85.00 \\ \hline
3 & 82.50 & 82.50 & 82.50 & 82.50 & 82.50 \\ \hline
4 & 81.67 & 83.33 & 82.50 & 85.83 & 84.17 \\ \hline
5 & 76.67 & 77.50 & 77.50 & 79.17 & 79.17 \\ \hline
Average & 81.83 & 82.33 & 82.17 & 82.83 & 82.67 \\ \hline
\end{tabular}
\end{center}
\end{table}
\begin{table}
\caption{5-sample-set} \label{5-sample-set}
\footnotesize
\begin{center}
\begin{tabular}{|c|r|r|r|r|r|} \hline
\multicolumn{1}{|p{40pt}|}{Learning Data Set} & \multicolumn{1}{p{25pt}|}{$\sh^{U}$} & \multicolumn{1}{p{25pt}|}{$\sh^{SDS}$} & \multicolumn{1}{p{25pt}|}{$\sh^{KG}$} & \multicolumn{1}{p{25pt}|}{$\sh^{MA\_1}$} & \multicolumn{1}{p{25pt}|}{$\sh^{MA\_2}$}\\ \hline
1 & 66.67 & 71.85 & 68.89 & 75.56 & 75.56 \\ \hline
2 & 78.52 & 80.00 & 78.52 & 85.19 & 82.96 \\ \hline
3 & 41.48 & 41.48 & 41.48 & 44.44 & 42.96 \\ \hline
4 & 43.70 & 46.67 & 45.93 & 53.33 & 50.37 \\ \hline
5 & 88.89 & 88.15 & 88.89 & 92.59 & 90.37 \\ \hline
6 & 73.33 & 78.52 & 77.78 & 89.63 & 88.15 \\ \hline
7 & 64.44 & 68.89 & 67.41 & 73.33 & 71.85 \\ \hline
8 & 73.33 & 75.56 & 72.59 & 82.96 & 79.26 \\ \hline
9 & 73.33 & 75.56 & 72.59 & 82.96 & 79.26 \\ \hline
10 & 69.63 & 72.59 & 71.85 & 82.22 & 77.78 \\ \hline
Average & 67.33 & 69.93 & 68.59 & 76.22 & 73.85 \\ \hline
\end{tabular}
\end{center}
\end{table}
We observe that 1) in each group, the largest eigenvalue are about 6 times as large as the second largest eigenvalue, 2) the second largest eigenvalue is about 3--7 times as large as the smallest eigenvalue. We are interested in the performance of $\s^{MA_\lossd},\ \lossd=1,2$, with the population eigenvalues in (\ref{iris_data_eigenvalues}) which are considered as a deviation from $(\infty, c, c, c)$, the ideal eigenvalues for $\s^{MA_i},\ i=1,2$. 

We made three types of cross validations.
\begin{enumerate}
\item Leave-one-out: For a chosen $(i,j),\ i=1,2,3,\ j=1,\ldots,50$, leave $\bm{x}_j^{(i)}$ out from the whole data to be a test data, and use the rest as a learning data set. We repeat this trial for every possible $(i,j)$. Consequently 150 trials were carried out.
\item 10-sample-set: First choose $\bm{x}_1^{(i)},\ldots,\bm{x}_{10}^{(i)},\ i=1,2,3$, as a learning data set and use all the rest as a test data. Next use $\bm{x}_{11}^{(i)},\ldots,\bm{x}_{20}^{(i)},\ i=1,2,3$, as a learning data set and the others as a test data. Repeatedly change a learning data set until every data is used once as a learning data. Totally we carried out $600(=120\times 5)$ trials.
\item 5-sample-set: First choose $\bm{x}_1^{(i)},\ldots,\bm{x}_5^{(i)},\ i=1,2,3$, as a learning data set and use all the rest as a test data. Next use $\bm{x}_6^{(i)},\ldots,\bm{x}_{10}^{(i)},\ i=1,2,3$, as a learning data set and the others as a test data. Repeatedly change a learning data set until every data is used once as a learning data. Totally we carried out $1350(=135\times 10)$ trials.
\end{enumerate}

We summarize the result on the correct classification percentage (``C.C.P.'' for abbreviation) of each discriminant function.
\begin{enumerate}
\item Leave-one-out: All the discriminant functions returned the same classification for every test data and scored 96.67\% of C.C.P. The misclassification occurred at the sample $\bm{x}^{(2)}_{19}$, $\bm{x}^{(2)}_{21}$, $\bm{x}^{(2)}_{23}$, $\bm{x}^{(2)}_{34}$, $\bm{x}^{(3)}_{32}$. With as much as 49 learning data, all the discrininant functions work quite correctly and make no differences among the functions.

\item 10-sample-set: See Table \ref{10-sample-set} for the C.C.P. in
  each learning data set and the average. Depending on the learning
  data set, different discriminant functions records the best
  C.C.P, but the margins are small and negligible. 
  It seems that even 10-sample-learning set is too large to
  differentiate the functions.
 
\item 5-sample-set: See Table \ref{5-sample-set} for the C.C.P. in
  each learning data set and the average. In every learning data set,
  the functions based on $\sh^{MA_\lossd},\ \lossd=1,2$, outperform the other
  functions. Especially $\sh^{MA_1}$ always keeps the highest C.C.P.
  In total, $\sh^{MA_1}$ and $\sh^{MA_2}$ record better C.C.P. than
  $\sh^U$ by 8.89\% and 6.52\% respectively, while the margins of
  $\sh^{SDS}$ and $\sh^{KG}$ over $\sh^{U}$ are respectively 2.60\% and
  1.26\%.
\end{enumerate}

\appendix
\section{Appendix}
\label{sec:appendix}
\subsection{Proof of Lemma 1}
\label{subsec:proof_lemma1}
In the following, $c_i\ (i=1,\ldots,7)$ represents some constant independent of $\aa,\bb$. 

The random variables $\lll=(l_1,\ldots,l_p)$ and $\gtilde=\ggammat\g$ have the following joint 
density function with respect to the product measure between Lebesgue measure 
on ${\cal L}$ and the invariant probability $\mu$ on $\opp$.
$$
 c_1\:\prodp \lambda_i^{-\frac{n}{2}} 
 \lff\:\lfs\:\etr \left(-\frac{1}{2}\gtilde\Lll\gtildet\Llambda^{-1}\right).
$$
We have 
\begin{eqnarray*}
E[x(\g,\lll,\llambda,\aa,\bb)]
&=&E[x(\ggamma\gtilde,\lll,\llambda,\aa,\bb)]\nonumber\\
&=&c_1\prodp \lambda_i^{-\frac{n}{2}}
\int_{\cal L}\int_{\op^+}x(\ggamma\g,\lll,\llambda,\aa,\bb)\lff\:\lfs\\
& &\hspace{3cm}\times   \etr \left(-\frac{1}{2}\g\Lll\ggt\Llambda^{-1}\right)d\mu(\g)d\lll.
\end{eqnarray*}
Using the finite open cover $\os$, $\covernum=0,\ldots,\covertotal$,
in Subsection \ref{subsec:localcoord}, we have
\begin{equation}
\label{comp_exp_x}
E[x(\g,\lll,\llambda,\aa,\bb)]=\sum_{\covernum=0}^\covertotal I_\covernum,
\end{equation}
where
\begin{eqnarray*}
I_\covernum&=&c_1\prodp \lambda_i^{-\frac{n}{2}}
\int_{\cal L}\int_{\opp} \iota_\covernum (\g) x(\ggamma\g,\lll,\llambda,\aa,\bb)\lff\:\lfs\\
& &\hspace{3cm}\times   \etr \left(-\frac{1}{2}\g\Lll\ggt\Llambda^{-1}\right)d\mu(\g)d\lll,
\end{eqnarray*}
First we consider $I_0$. Let $M$ denote the support of $\iota_0$. From (\ref{bound_cond_x}),
\begin{eqnarray}
 |I_0|&\leq&c_1\:\dts \int_{\cal L}\int_{M}\:|x(\ggamma\g,\lll,\llambda,\aa,\bb)|\:\lff\:\lfs\nonumber\\
& & \hspace{3cm}\times\etr
  \left(-\frac{1}{2}\g\Lll\ggt\Llambda^{-1}\right)d\mu(\g)d\lll 
  \nonumber \\
&\le& c_1\:b\:\dts \int_{\cal L}\int_{M}\:\lff\:\lfs\:
  \etr
  \left(-\frac{1}{2}\g\Lll\ggt\bar{\Llambda}^{-1}\right)d\mu(\g)d\lll 
  \nonumber \\
&=&c_2\:P\left(\gtilde \in M\:|\: \s=\ggamma\bar{\Llambda}
\ggammat\right),
\end{eqnarray}
where $\bar{\Llambda}=(1-2a)^{-1}\Llambda.$ Note $\opp\setminus M$ is an open set including ${\cal O}(m,p-m)$, hence 
by 2 of Theorem \ref{thm:asympt_g21}, $\lim_{\bdiva\rightarrow 0}P\left(\gtilde \in \opp\setminus M\:|\: \s=\ggamma\bar{\Llambda}
\ggammat\right)=1$, which means
$$P\left(\gtilde \in M\:|\: \s=\ggamma\bar{\Llambda}
\ggammat\right)\rightarrow 0 $$
as $\bdiva\rightarrow 0$. Therefore 
\begin{equation}
\label{vanish_I0}
\lim_{\bdiva\rightarrow 0}I_0=0.
\end{equation}

Now we focus ourselves on $I_\covernum$, 
$\covernum=1,\ldots,\covertotal$. Since $\mu$ is invariant and the support of $\iota_\covernum(\g)$ is contained in $O^{(\covernum)}$, we have
\begin{eqnarray*}
I_\covernum&=&c_1\prodp \lambda_i^{-\frac{n}{2}}
\int_{\cal L}\int_{V}
\iota_\covernum(\hh^{(\covernum)}\g)\;
x(\ggamma\hh^{(\covernum)}\g,\lll,\llambda,\aa,\bb)\lff\:\lfs\\
& &\hspace{4cm}\times   \etr \left(-\frac{1}{2}\hh^{(\covernum)}\g\Lll\ggt
\hh^{(\covernum)}{}'\Llambda^{-1}\right)d\mu(\g)d\lll.
\end{eqnarray*}
We want to express the integral with respect to $d\mu(\g)$ in terms of
the local coordinates $\uu$ on $U$. 
It is well known that the invariant measure 
$d\mu(\g)$ has the exterior differential form expression
\begin{equation}
\label{eq:invariant-measure-in-differentail-form}
c_3\bigwedge_{i>j}\ggg_j' d\ggg_i,
\end{equation}
where $\ggg_i$ is the $i$th column of $\g$.  Substituting the differential
\begin{eqnarray*}
d g_{ij} &=& d u_{ij}, \quad i > j, \\
d g_{ij} &=& \sum_{k > l} \frac{\partial g_{ij}}{\partial u_{kl}} d
  u_{kl}, \quad  i \le j,
\end{eqnarray*}
into (\ref{eq:invariant-measure-in-differentail-form}) and taking the
wedge product of the terms, we see that  
$$
\bigwedge_{i>j}\ggg_j' d\ggg_i = \pm J^*(\uu) \bigwedge_{i>j} du_{ij}, 
$$
where $J^*(\uu)$ is the Jacobian expressing the Radon-Nikodym derivative
of the measure on $U$ induced from the invariant measure on
$\opp$ with respect to the Lebesgue measure on $R^{\frac{p(p-1)}{2}}$.
An explicit form of $J^*(\uu)$ for small dimension $p$ is discussed in
Appendix B in Takemura and Sheena (2005). 
Since $J^*(\uu)$ is a $C^\infty$ function on $\bar{U}$, it is
bounded and has a finite limit as $\uu\rightarrow \zz$.
By the change of variables $(\lll,\g)\rightarrow (\lll,\uu)$, 
$I_\covernum$ is written as
\begin{eqnarray*}
I_\covernum&=&c_4\prodp \lambda_i^{-\frac{n}{2}}
\int_{\cal L}\int_{U}
\iota_\covernum(\hh^{(\covernum)}\g(\uu))\;
x(\ggamma\hh^{(\covernum)}\g(\uu),\lll,\llambda,\aa,\bb)\lff\:\lfs\\
& &\hspace{3cm}\times   \etr
\left(-\frac{1}{2}\hh^{(\covernum)}\g(\uu)\Lll\g'(\uu)
\hh^{(\covernum)}{}'\Llambda^{-1}\right)J^*(\uu)d\uu d\lll,
\end{eqnarray*}
Consider further coordinate transformation $(\lll,\uu)\rightarrow
(\ddd,\qq)$ for each $\covernum$.  Notice
\begin{equation}
\label{trans1}
\prod_{i=1}^pl_i^{\frac{n-p-1}{2}}=\Bigl(\prod_{i=1}^p
d_j^{\frac{n-p-1}{2}}\Bigr)\aa^{\frac{m(n-p-1)}{2}}\bb^{\frac{(p-m)(n-p-1)}{2}}, 
\end{equation}
\begin{eqnarray}
\label{trans2}
\prod_{j<i}(l_j-l_i)
&=&\aa^{\frac{m(m-1)}{2}}\bb^{\frac{(p-m)(p-m-1)}{2}}\prod_{j\leq m <i}(\aa d_j-\bb d_i)
\prod_{j<i\leq m}(d_j-d_i)
\prod_{m<j<i}(d_j-d_i) \nonumber\\
&=&\prod_{j\leq m <i}\Bigl(1-\frac{\bb d_i}{\aa
  d_j}\Bigr)\prod_{j<i\leq m}(d_j-d_i)\prod_{m<j<i}(d_j-d_i) \prod_{j=1}^m
d_j^{\:p-m}\nonumber\\
& &\qquad \times
\aa^{m(p-m)+\frac{m(m-1)}{2}}\bb^{\frac{(p-m)(p-m-1)}{2}}, 
\end{eqnarray} 
and
\begin{eqnarray}
\label{trans3}
\lefteqn{
\tr\hh^{(\covernum)}\g(\uu)\Lll\g'(\uu)\hh^{(\covernum)}{}'\Llambdain}\nonumber\\
&=&\tr\left\{
\left(
\begin{array}{cc}
\hh_1^{(\covernum)}\g_{11}(\uu) & \hh_1^{(\covernum)}\g_{12}(\uu) \\
\hh_2^{(\covernum)}\g_{21}(\uu)& \hh_2^{(\covernum)}\g_{22}(\uu) 
\end{array}
\right)
\mbox{diag}(l_1,\ldots,l_p)
\right.
\nonumber\\
&&\qquad \times 
\left.
\left(
\begin{array}{cc}
\g'_{11}(\uu)\hh_1^{(\covernum)}{}' & \g'_{21}(\uu)\hh_2^{(\covernum)}{}' \\
\g'_{12}(\uu)\hh_1^{(\covernum)}{}'& \g'_{22}(\uu)\hh_2^{(\covernum)}{}' 
\end{array}
\right)
\mbox{diag}(\lambda_1^{-1},\ldots,\lambda_p^{-1})
\right\}
\nonumber\\
&=& \tr
\hh_1^{(\covernum)}\g_{11}(\uu)\dd_1\g'_{11}(\uu)\hh_1^{(\covernum)}{}'\xxxi_1^{-1}
+\tr \hh_2^{(\covernum)}\g_{22}(\uu)\dd_2\g'_{22}(\uu)\hh_2^{(\covernum)}{}'\xxxi_2^{-1}
\nonumber\\
& &\ + \tr\bm{Q}_{21}\bm{Q}'_{21} +\aa^{-1}\bb \tr \hh_1^{(\covernum)}\g_{12}(\uu)\dd_2\g'_{12}(\uu)\hh_1^{(\covernum)}{}'\xxxi_1^{-1},
\end{eqnarray}
where $\uu$ is actually the abbreviation for
$\uu(\ddd,\qq,\xxi,\aa,\bb)$
defined by (\ref{def_u( )}). 
For notational simplicity we use the same abbreviation
$\uu=\uu(\ddd,\qq,\xxi,\aa,\bb)$ for the rest of this proof.
From (\ref{jac_lu_tilq}), (\ref{exp_x}), (\ref{trans1}), (\ref{trans2}) and (\ref{trans3}), we have
$$
I_\covernum
=c_5\int_{R^{p(p-1)/2}}\int_{R^p_+}
\iota_\covernum(\hh^{(\covernum)}
\g(\uu))x(\ddd,\qq,\xxi,\aa,\bb;\ggamma,\hh^{(\covernum)}) h(\ddd,\qq,\xxi,\aa,\bb)\:d\ddd d\qq
$$
where $R_+^p=\{\ddd\:|\:d_i>0,\ i=1,\ldots,p\}$ and $h(\ddd,\qq,\xxi,\aa,\bb)$ is defined as follows;
\begin{eqnarray*}
h(\ddd,\qq,\xxi,\aa,\bb)
&=&I(\uu\in U) J^*(\uu) I(\ddd_1 \in {\cal D}_1,\ \ddd_2 \in {\cal D}_2,
\  (\ddd_1,\ddd_2)\in  {\cal D}_3)\\
& &\ \times\prod_{i=1}^m d_i^{\frac{n-m-1}{2}} \prod_{i=m+1}^p d_i^{\frac{n-p-1}{2}}
\prod_{j<i\leq m}(d_j-d_i) \prod_{m<j<i}(d_j-d_i)\prod_{j\leq m <i}\Bigl(1-\frac{\bb d_i}{\aa d_j}\Bigr)\\
 & &\
 \times\exp\Bigl(-\frac{1}{2}\tr\sum_{s=1}^2 \hh_s^{(\covernum)}\g_{ss}(\uu)\dd_s\g'_{ss}(\uu)
 \hh_s^{(\covernum)}{}'\xxxi_s^{-1}\Bigr)\\
& &\ \times\etr\Bigl(-\frac{1}{2}\bm{Q}_{21}\bm{Q}'_{21}\Bigr)
\times\etr\Bigl(-\frac{\bb}{2\aa}\hh_1^{(\covernum)}\g_{12}(\uu)\dd_2\g'_{12}(\uu)
\hh_1^{(\covernum)}{}'\xxxi_1^{-1}\Bigr).
\end{eqnarray*}

We will show that 
$$
 \iota_\covernum(\hh^{(\covernum)}\g(\uu))x(\ddd,\qq,\xxi,\aa,\bb;\ggamma,\hh^{(\covernum)})h(\ddd,\qq,\xxi,\aa,\bb)
$$
is bounded in $(\aa,\bb)$. First $I(\uu\in U) J^*(\uu)\leq K$
for some $K\,(>0)$ since $J^*(\uu)$ is bounded on the compact set $\bar{U}$. Clearly 
$$
0\leq I(\ddd_1 \in {\cal D}_1,\ \ddd_2 \in {\cal D}_2,
\  (\ddd_1,\ddd_2)\in  {\cal D}_3)\prod_{j\leq m <i}
\Bigl(1-\frac{\bb d_i}{\aa d_j}\Bigr)\leq 1.
$$
From the condition (\ref{bound_cond_x}), we have
\begin{eqnarray*}
|x(\ddd,\qq,\xxi,\aa,\bb;\ggamma,\hh^{(\covernum)})| 
&=&|x(\ggamma\hh^{(\covernum)}\g(\uu),\lll,\llambda,\aa,\bb)|\\
&\leq & b\: \etr(a\hh^{(\covernum)}\g(\uu)\Lll\g'(\uu)\hh^{(\covernum)}{}'\Llambdain)\ a.e. \mbox{ in } (\ddd,\qq).
\end{eqnarray*}
Therefore 
\begin{eqnarray}
\lefteqn{| \iota_\covernum(\hh^{(\covernum)}\g(\uu))x(\ddd,\qq,\xxi,\aa,\bb;\ggamma,\hh^{(\covernum)})h(\ddd,\qq,\xxi,\aa,\bb)|}\nonumber\\
&\leq&
c_6\:I(\uu\in U)\prod_{i=1}^m d_i^{\:\frac{n-m-1}{2}}
\prod_{i=m+1}^p d_i^{\:\frac{n-p-1}{2}}
\prod_{j<i\leq m}|d_j-d_i| \prod_{m<j<i}|d_j-d_i| 
\nonumber\\
& &\ \times\exp\Bigl(-\frac{1-2a}{2}\tr\sum_{s=1}^2 
\hh_s^{(\covernum)}\g_{ss}(\uu)\dd_s\g'_{ss}(\uu)\hh_s^{(\covernum)}{}'\xxxi_s^{-1}\Bigr)
\nonumber \\
& &\ 
\times\etr\Bigl(-\frac{1-2a}{2}\bm{Q}_{21}\bm{Q}'_{21}\Bigr).
\label{upbound_xh}
\end{eqnarray}
Note that 
$$
I(\uu \in U)
\leq I(\uu\in C_\epsilon)
\leq I(|u_{ij}|=|q_{ij}|<\epsilon,\ 1\leq j<i \leq m,\ m<j<i\leq p).
$$
Choose some $\bar{\xi}$ such that $\bar{\xi}>\xi_i,\ i=1,\ldots,p$.
Consequently the left-hand side of (\ref{upbound_xh}) is bounded by $\bar{h}(\ddd,\qq)$, where
\begin{eqnarray*}
\bar{h}(\ddd,\qq)
&=& c_6\:I(|q_{ij}|<\epsilon,\ 1\leq j<i \leq m,\ m<j<i\leq p)\\
& &\ \times \prod_{i=1}^m d_i^{\:\frac{n-m-1}{2}} \prod_{i=m+1}^p d_i^{\:\frac{n-p-1}{2}}\prod_{j<i\leq m}|d_j-d_i| \prod_{m<j<i}|d_j-d_i| \\
& &\ \times \exp\Bigl(-\frac{\bar{\xi}^{-1}}{2}(1-2a)\sum_{i=1}^pd_i\Bigr)
\times \etr\Bigl(-\frac{1-2a}{2}\bm{Q}_{21}\bm{Q}'_{21}\Bigr).
\end{eqnarray*}
Let $\nu_1=m(m-1)/2,\ \nu_2=(p-m)(p-m-1)/2,\  \nu_3=m(p-m)$. We have
\begin{eqnarray*}
\int_{R_+^p}\:\int_{R^{p(p-1)/2}}\bar{h}(\ddd,\qq)\:d\qq\:d\ddd
&=&
\int_{R_+^p}\int_{R^{\nu_3}}\int_{R^{\nu_2}}\int_{R^{\nu_1}}
\bar{h}(\ddd,\qq)\: d\qq_{11}\: d\qq_{22}\: d\qq_{21}\: d\ddd\\
&=& c_6 \int_{R_+^p}\prod_{i=1}^m d_i^{\:\frac{n-m-1}{2}} \prod_{i=m+1}^p d_i^{\:\frac{n-p-1}{2}}
\prod_{j<i\leq m}|d_j-d_i| \prod_{m<j<i}|d_j-d_i| \\
& & \times \exp\Bigl(-\frac{\bar{\xi}^{-1}}{2}(1-2a)\sum_{i=1}^pd_i\Bigr)d\ddd
\times \int_{R^{\nu_3}}
\etr\Bigl(-\frac{1-2a}{2}\bm{Q}_{21}\bm{Q}'_{21}\Bigr)\;d\qq_{21}\\
& &\times \int_{(-\epsilon,\epsilon)^{\nu_1}} 1 \:d\qq_{11}\
\int_{(-\epsilon,\epsilon)^{\nu_2}} 1 \:d\qq_{22} < \infty .
\end{eqnarray*}
The integrability of $\bar{h}(\ddd,\qq)$ guarantees the use of the dominated convergence theorem; From (\ref{asympt_G}) and (\ref{conv_x})
\begin{eqnarray*}
\lim_{\bdiva\rightarrow 0}I_\covernum
&=&c_5\int_{R^{p(p-1)/2}}\int_{R^p_+} \lim_{\bdiva\rightarrow 0}\iota_\covernum(\hh^{(\covernum)}\g(\uu))\lim_{\bdiva\rightarrow 0}x(\ddd,\qq,\xxi,\aa,\bb;\ggamma,\hh^{(\covernum)})\\
& &\hspace{3cm}\times \lim_{\bdiva\rightarrow 0}h(\ddd,\qq,\xxi,\aa,\bb)\:d\ddd d\qq\\
&=&c_5\int_{R^{p(p-1)/2}}\int_{R^p_+} \iota_\covernum(\hh^{(\covernum)}\g(\qq_{11},\qq_{22},\bm{0}))\:\bar{x}_{\ggamma}(\hh^{(\covernum)}\g(\qq_{11},\qq_{22},\bm{0}),\ddd,\qqq_{21},\xxi) \\
& &\hspace{3cm} \times\lim_{\bdiva\rightarrow 0}h(\ddd,\qq,\xxi,\aa,\bb)\:d\ddd d\qq.
\end{eqnarray*}
We consider $\lim_{\bdiva\rightarrow 0}h(\ddd,\qq,\xxi,\aa,\bb).$ First notice that 
\begin{eqnarray*}
&&\lim_{\bdiva\rightarrow 0}I(\ddd_1 \in {\cal D}_1,\ \ddd_2 \in {\cal D}_2,
\  (\ddd_1,\ddd_2)\in  {\cal D}_3)=I(\ddd_1\in {\cal D}_1)I(\ddd_2\in {\cal D}_2),\\
&&\lim_{\bdiva\rightarrow 0}
\prod_{j\leq m <i}\Bigl(1-\frac{\bb d_i}{\aa d_j}\Bigr)=1.
\end{eqnarray*}
From (\ref{asym_u}), we find
\begin{eqnarray*}
&&\lim_{\bdiva\rightarrow 0}J^*(\uu)=J^*(\qq_{11},\qq_{22},\bm{0}),\\
&&\lim_{\bdiva\rightarrow 0}I(\uu\in U)=I((\qq_{11},\qq_{22})=(\uu_{11},\uu_{22})\in U_0),
\end{eqnarray*}
where $U_0 = \{(\uu_{11},\uu_{22}) |(\uu_{11},\uu_{22},\bm{0})\in U\}$ 
denotes the slice of $U$ by $\uu_{12}=0,$ 
and that
\begin{eqnarray*}
&&\lim_{\bdiva\rightarrow 0}\g_{11}(\uu)=\g_{11}(\qq_{11},\qq_{22},\bm{0})\in {\cal O}^+(m),\\
&&\lim_{\bdiva\rightarrow 0}\g_{22}(\uu)=\g_{22}(\qq_{11},\qq_{22},\bm{0})\in {\cal O}^+(p-m),\\
&&\lim_{\bdiva\rightarrow 0}\g_{21}(\uu)=\bm{0},\\
&&\lim_{\bdiva\rightarrow 0}\etr\Bigl(-\frac{\beta}{2\alpha}\hh_1^{(\covernum)}\g_{12}(\uu)\dd_2\g'_{12}(\uu)\hh_1^{(\covernum)}{}'\xxxi_1^{-1}\Bigr)=1.
\end{eqnarray*}
Since $d\mu$ is invariant, especially w.r.t.\ both of the transformations
\begin{equation}
\label{trans_hh}
\g\rightarrow \diag(\hh_1,\hh_2)\g,\qquad \g\rightarrow \g\diag(\hh_1,\hh_2),
\end{equation}
the measure on $U_0$ given by
\begin{equation}
\label{indu_dist_U}
J^*(\qq_{11},\qq_{22},\bm{0}) d\qq_{11}d\qq_{22}
\end{equation}
induces the invariant measure on $V_0$, the slice of $V$
by $\g_{12}=0$, w.r.t.\ %
(\ref{trans_hh}) 
through 
\begin{equation}
\label{diag_g(u)}
\g(\qq_{11},\qq_{22},\bm{0})=
\diag(\g_{11}(\qq_{11},\qq_{22},\bm{0}), \g_{22}(\qq_{11},\qq_{22},\bm{0})).
\end{equation}
If $\g_{11}$ and $\g_{22}$ independently follow the invariant probability distributions respectively on ${\cal O}^+(m)$ and ${\cal O}^+(p-m)$, then the distribution on $V_0$ given by
\begin{equation}
\label{def_g0}
\g_0=
\diag(\g_{11}, \g_{22})
\end{equation}
is also invariant w.r.t.\ the transformations (\ref{trans_hh}), hence must be proportional to the above-mentioned distribution on $V_0$ given by (\ref{diag_g(u)}) and (\ref{indu_dist_U}). 
Consequently 
\begin{eqnarray*}
\lim_{\bdiva\rightarrow 0}I_\covernum
&=&c_6\int_{R^{m(p-m)}}\:\int_{R_+^p}\:\int_{V_0}\iota_\covernum(\hh^{(\covernum)}\g_0)\bar{x}_{\ggamma}(\hh^{(\covernum)}\g_0,\ddd,\qqq_{21},\xxi)
 I(\ddd_1 \in {\cal D}_1)I(\ddd_2 \in {\cal D}_2)\\
& &\hspace{2.5cm}\times\prod_{i=1}^m d_i^{\:\frac{n-m-1}{2}} \prod_{i=m+1}^p d_i^{\:\frac{n-p-1}{2}}\prod_{j<i\leq m}(d_j-d_i) \prod_{m<j<i}(d_j-d_i)\\
& &\hspace{2.5cm}\times
\exp\Bigl(-\frac{1}{2}\tr\sum_{s=1}^2 \hh_s^{(\covernum)}\g_{ss}\dd_s\g'_{ss}
\hh_s^{(\covernum)}{}'\xxxi_s^{-1}\Bigr)\\
& &\hspace{2.5cm}\times\etr\Bigl(-\frac{1}{2}\bm{Q}_{21}\bm{Q}'_{21}\Bigr) d\mu_1(\g_{11})\:d\mu_2(\g_{22})
\:d\ddd\:d\qq_{21},
\end{eqnarray*}
where $\g_0$ is given by (\ref{def_g0}), and $\mu_1$, $\mu_2$ are the invariant probability measures respectively on 
${\cal O}^+(m)$ and ${\cal O}^+(p-m).$ 

Let $O_0^{(\covernum)}$ denote the slice of $O^{(\covernum)}$ by $\g_{12}=\bm{0}.$ Since 
$O^{(\covernum)}=\hh^{(\covernum)}V$, $O_0^{(\covernum)}=\hh^{(\covernum)}V_0$.
Consequently for each $1\leq \covernum \leq \covertotal$,
\begin{eqnarray*}
\lim_{\bdiva\rightarrow 0}I_\covernum
&=&c_6\int_{R^{m(p-m)}}\:\int_{R_+^p}\:\int_{O_0^{(\covernum)}}\iota_\covernum(\g_0)\;\bar{x}_{\ggamma}(\g_0,\ddd,\qqq_{21},\xxi)
I(\ddd_1 \in {\cal D}_1)I(\ddd_2 \in {\cal D}_2)\nonumber\\
& &\times\prod_{i=1}^m d_i^{\:\frac{n-m-1}{2}}\prod_{i=m+1}^p d_i^{\:\frac{n-p-1}{2}}\prod_{j<i\leq m}(d_j-d_i) \prod_{m<j<i}(d_j-d_i)\nonumber\\
& &\times\exp\Bigl(-\frac{1}{2}\tr\sum_{s=1}^2
\g_{ss}\dd_s\g'_{ss}\xxxi_s^{-1}\Bigr)
\etr\Bigl(-\frac{1}{2}\bm{Q}_{21}\bm{Q}'_{21}\Bigr)
d\mu_1(\g_{11})\:d\mu_2(\g_{22})
\:d\ddd\:d\qq_{21}.
\end{eqnarray*}
Note that $\bigcup_{\covernum=1}^\covertotal 
 O^{(\covernum)}_0={\cal O}(m,p-m)$ and $\iota_\covernum(\g_0)$ vanishes on ${\cal O}(m,p-m)\setminus O^{(\covernum)}_0$. Therefore we have 
\begin{eqnarray}
\lim_{\bdiva\rightarrow 0}I_\covernum
&=&c_6\int_{R^{m(p-m)}}\:\int_{R_+^p}\:\int_{{\cal O}(m,p-m)}\iota_\covernum(\g_0)\bar{x}_{\ggamma}(\g_0,\ddd,\qqq_{21},\xxi)
I(\ddd_1 \in {\cal D}_1)I(\ddd_2 \in {\cal D}_2)\nonumber\\
& &
\times\prod_{i=1}^m
d_i^{\:\frac{n-m-1}{2}}\prod_{i=m+1}^p
d_i^{\:\frac{n-p-1}{2}}\prod_{j<i\leq m}(d_j-d_i)
\prod_{m<j<i}(d_j-d_i)
\nonumber \\
& &\times\exp\Bigl(-\frac{1}{2}\tr\sum_{s=1}^2 \g_{ss}\dd_s\g'_{ss}\xxxi_s^{-1}\Bigr)
\etr\Bigl(-\frac{1}{2}\bm{Q}_{21}\bm{Q}'_{21}\Bigr)
d\mu_1(\g_{11})\:d\mu_2(\g_{22})
\:d\ddd\:d\qq_{21}.
\label{lim_Is}
\end{eqnarray}

From (\ref{comp_exp_x}), (\ref{vanish_I0}) and ({\ref{lim_Is}), we have
\begin{eqnarray}
\label{lim_E[x(g,l...)]}
\lefteqn{\lim_{\bdiva\rightarrow 0} E[x(\g,\lll,\llambda,\aa,\bb)]}\nonumber\\
&=&c_6\int_{R^{m(p-m)}}\:\int_{R_+^p}\:\int_{{{\cal O}^+(p-m)}}\:\int_{{{\cal O}^+(m)}}\bar{x}_{\ggamma}(\g_0,\ddd,\qqq_{21},\xxi) I(\ddd_1 \in {\cal D}_1)I(\ddd_2 \in {\cal D}_2)
\nonumber\\
& &\times\prod_{i=1}^m d_i^{\:\frac{n-m-1}{2}} \prod_{i=m+1}^p d_i^{\:\frac{n-p-1}{2}}\prod_{j<i\leq m}(d_j-d_i) \prod_{m<j<i}(d_j-d_i)\\
&&
\times\exp\Bigl(-\frac{1}{2}\tr\sum_{s=1}^2
 \g_{ss}\dd_s\g'_{ss}\xxxi_s^{-1}\Bigr)
\times\etr\Bigl(-\frac{1}{2}\bm{Q}_{21}\bm{Q}'_{21}\Bigr) d\mu_1(\g_{11})\:d\mu_2(\g_{22})
\:d\ddd\:d\qq_{21}\nonumber.
\end{eqnarray}
Under the distribution (\ref{asympt_dist_lemma1}) and the spectral
decompositions (\ref{decomp_W_ii}), the joint density function of
$(\ddd_{1}, \g_{11})$ ($(\ddd_{2}, \g_{22})$) with respect to the
product measure of Lebesgue measure on $R^m_+$ ($R^{p-m}_+$) and the
invariant probability measure $\mu_1$ ($\mu_2$) on 
${\cal O}^+(m)$ (${\cal O}^+({p-m})$) is given by the following
functions, $F_1(\ddd_{1}, \g_{11})$ ($F_2(\ddd_{2}, \g_{22})$);
\begin{eqnarray*}
F_1(\g_{11},\ddd_1)&=&K_1|\xxxi_1|^{-\frac{n}{2}}\prod_{i=1}^md_i^{\:\frac{n-m-1}{2}}\prod_{1\leq j<i \leq m}(d_j-d_i)
\etr\Bigl(
-\frac{1}{2}\g_{11}\dd_1\g'_{11}\xxxi_1^{-1}\Bigr)\\
F_2(\g_{22},\ddd_2)&=&K_2|\xxxi_2|^{-\frac{n-m}{2}}\prod_{i=m+1}^pd_i^{\:\frac{n-p-1}{2}}\prod_{m<j<i\leq p}(d_j-d_i)
\etr\Bigl(
-\frac{1}{2}\g_{22}\dd_2\g'_{22}\xxxi_2^{-1}\Bigr),
\end{eqnarray*}
with $K_1, K_2$ as normalizing constants.
The density function of $\bm{Z}_{21}$ is given by 
$$
F_3(\bm{z}_{21})=K_3 \etr(-\frac{1}{2}\bm{Z}_{21}\bm{Z}'_{21}),
$$
where $K_3$ is a normalizing constant. Using $F_1(\g_{11},\ddd_1), F_2(\g_{22},\ddd_2), F_3(\bm{z}_{21})$, we can rewrite the right-hand side of (\ref{lim_E[x(g,l...)]}) as
\begin{eqnarray*}
&&c_7\int_{R^{m(p-m)}}\:\int_{{\cal D}_2}\:\int_{{\cal D}_1}\:\int_{{{\cal O}^+(p-m)}}\:\int_{{{\cal O}^+(m)}}\bar{x}_{\ggamma}(\g_0,(\ddd_1,\ddd_2),\bm{Z}_{21},\xxi)\\
&&\hspace{2cm}\times F_1(\g_{11},\ddd_1)F_2(\g_{22},\ddd_2)F_3(\bm{z}_{21})
d\mu_1(\g_{11})\:d\mu_2(\g_{22})
\:d\ddd_1\:d\ddd_2\:d\bm{z}_{21}.
\end{eqnarray*}
If we consider the special case $x(\g,\lll,\llambda,\aa,\bb)=1$, we notice that $c_7=1.$
\subsection{Proof of Lemma 2}
\label{subsec:proof_lemma2}
Using Lemma 1, we will calculate 
\begin{eqnarray*}
&&\lim_{\bdiva \rightarrow 0}E[\tr\g\Lllh\cc\Lllh\ggt\ggamma\Llambdain\ggammat],\\
&&\lim_{\bdiva \rightarrow 0}E[\tr(\g\Lllh\cc\Lllh\ggt\ggamma\Llambdain\ggammat)^2].
\end{eqnarray*}
Let 
\begin{eqnarray*}
x_1(\g,\lll,\llambda,\aa,\bb)&=&\tr(\g\Lllh\cc\Lllh\ggt\ggamma\Llambdain\ggammat),\\
x_2(\g,\lll,\llambda,\aa,\bb)&=&\tr(\g\Lllh\cc\Lllh\ggt\ggamma\Llambdain\ggammat)^2,
\end{eqnarray*}
then
\begin{eqnarray*}
x_1(\ggamma\g,\lll,\llambda,\aa,\bb)
&=&\sumi \sumj \lambda_i^{-1} l_j c_j g_{ij}^2 \\
&\leq & (\max_j c_j)\sumi \sumj \lambda_i^{-1}l_j g_{ij}^2
=3(\max_j c_j) \tr\Bigl(\frac{1}{3}\g\Lll\ggt\Llambdain\Bigr)\\
&\leq& 3(\max_j c_j) \etr\Bigl(\frac{1}{3}\g\Lll\ggt\Llambdain\Bigr),\\
x_2(\ggamma\g,\lll,\llambda,\aa,\bb)
&=&\tr(\Llambdaih\g\Lllh\cc\Lllh\ggt\Llambdaih)^2\\
 &\leq&\Bigl(\tr\Llambdaih\g\Lllh\cc\Lllh\ggt\Llambdaih \Bigr)^2
 =\Bigl(\sumi \sumj \lambda_i^{-1} l_j c_j g_{ij}^2 \Bigr)^2\\
 &\leq & \Bigl((\max_j c_j)\sumi \sumj \lambda_i^{-1}l_j g_{ij}^2\Bigr)^2
 = \Bigl\{6(\max_j c_j) \tr\Bigl(\frac{1}{6}\g\Lll\ggt\Llambdain\Bigr)\Bigr\}^2\\
&\leq& \Bigl\{6(\max_j c_j) \etr\Bigl(\frac{1}{6}\g\Lll\ggt\Llambdain\Bigr)\Bigr\}^2 
=36(\max_j c_j)^2\etr\Bigl(\frac{1}{3}\g\Lll\ggt\Llambdain\Bigr),                 
\end{eqnarray*}                        
hence (\ref{bound_cond_x}) is satisfied for both $x_1$ and $x_2$.                  
Now let
$$
\bm{B}(\g,\lll,\llambda)=\Llambdaih\hh^{(\covernum)}\g\Lllh
$$
for each $\covernum$.  Then we have
\begin{eqnarray*}
x_1(\ggamma\hh^{(\covernum)}\g,\lll,\llambda,\aa,\bb)&=&\tr\bm{B}\cc\bm{B}',\\
x_2(\ggamma\hh^{(\covernum)}\g,\lll,\llambda,\aa,\bb)&=&\tr(\bm{B}\cc\bm{B}')^2.
\end{eqnarray*}
We notice that
\begin{eqnarray*}
\bm{B}&=&\Llambdaih\hh^{(\covernum)}\g(\uu)\Lllh\\
 &=&\left(
\begin{array}{cc}
\Llambdaih_1\hh_1^{(\covernum)}\g_{11}(\uu)\Lllh_1 & \Llambdaih_1\hh_1^{(\covernum)}\g_{12}(\uu)\Lllh_2\\
\Llambdaih_2\hh_2^{(\covernum)}\uuu_{21}\Lllh_1    & \Llambdaih_2\hh_2^{(\covernum)}\g_{22}(\uu)\Lllh_2
\end{array}
\right)\\
&=&\left(
\begin{array}{cc}
\xxxi_1^{-1/2}\hh_1^{(\covernum)}\g_{11}(\uu)\dd_1^{1/2} & \aa^{-1/2}\bb^{1/2}\xxxi_1^{-1/2}\hh_1^{(\covernum)}\g_{12}(\uu)\dd_2^{1/2}\\
\aa^{1/2}\bb^{-1/2}\xxxi_2^{-1/2}\hh_2^{(\covernum)}\uuu_{21}\dd_1^{1/2}    & \xxxi_2^{-1/2}\hh_2^{(\covernum)}\g_{22}(\uu)\dd_2^{1/2}
\end{array}
\right)\\
&=&\left(
\begin{array}{cc}
\xxxi_1^{-1/2}\hh_1^{(\covernum)}\g_{11}(\uu)\dd_1^{1/2} & \aa^{-1/2}\bb^{1/2}\xxxi_1^{-1/2}\hh_1^{(\covernum)}\g_{12}(\uu)\dd_2^{1/2}\\
\qqq_{21}   & \xxxi_2^{-1/2}\hh_2^{(\covernum)}\g_{22}(\uu)\dd_2^{1/2}
\end{array}
\right).
\end{eqnarray*}
Substitute $\uu(\ddd,\qq,\xxi,\aa,\bb)$ with $\uu$ in the last matrix
and denote it by  
$\bm{B}(\ddd,\qq,\xxi,\aa,\bb)$. 
Then 
\begin{eqnarray*}
x_1(\ddd,\qq,\xxi,\aa,\bb;\ggamma,\hh^{(\covernum)})
&=&\tr\bm{B}(\ddd,\qq,\xxi,\aa,\bb)\cc 
\bm{B}'(\ddd,\qq,\xxi,\aa,\bb),\\
x_2(\ddd,\qq,\xxi,\aa,\bb;\ggamma,\hh^{(\covernum)})
&=&\tr(\bm{B}(\ddd,\qq,\xxi,\aa,\bb)\cc 
\bm{B}'(\ddd,\qq,\xxi,\aa,\bb))^2.
\end{eqnarray*}
Therefore
\begin{eqnarray*}
\lim_{\bdiva\rightarrow 0}x_1(\ddd,\qq,\xxi,\aa,\bb;\ggamma,\hh^{(\covernum)})
&=&\tr\bar{\bm{B}}\cc \bar{\bm{B}}',\\
\lim_{\bdiva\rightarrow 0}x_2(\ddd,\qq,\xxi,\aa,\bb;\ggamma,\hh^{(\covernum)})
&=&\tr(\bar{\bm{B}}\cc \bar{\bm{B}}')^2,
\end{eqnarray*}
where
$$
\bar{\bm{B}}=\left(
\begin{array}{cc}
\bar{\bm{B}}_{11} & \bar{\bm{B}}_{12} \\
\bar{\bm{B}}_{21} & \bar{\bm{B}}_{22}
\end{array}
\right)=\lim_{\bdiva\rightarrow 0}\bm{B}(\ddd,\qq,\xxi,\aa,\bb)
$$
is given by
\begin{eqnarray*}
\bar{\bm{B}}_{11}&=&\xxxi_1^{-1/2}\hh_1^{(\covernum)}\g_{11}(\qq_{11},\qq_{22},\bm{0})\dd_1^{1/2}=
\xxxi_1^{-1/2}(\hh^{(\covernum)}\g(\qq_{11},\qq_{22},\bm{0}))_{11}\dd_1^{1/2},\\
\bar{\bm{B}}_{12}&=&\bm{0},\\
\bar{\bm{B}}_{21}&=&\qqq_{21},\\
\bar{\bm{B}}_{22}&=&\xxxi_2^{-1/2}\hh_2^{(\covernum)}\g_{22}(\qq_{11},\qq_{22},\bm{0})\dd_2^{1/2}=
\xxxi_2^{-1/2}(\hh^{(\covernum)}\g(\qq_{11},\qq_{22},\bm{0}))_{22}\dd_2^{1/2},
\end{eqnarray*}
because of (\ref{asympt_G}).
By straightforward calculation we have
\begin{eqnarray*}
\tr(\bar{\bm{B}}\cc\bar{\bm{B}}')
&=&\tr\bar{\bm{B}}_{11}\cc_1\bar{\bm{B}}_{11}'+\tr\bar{\bm{B}}_{12}\cc_2
\bar{\bm{B}}_{12}'
+\tr\bar{\bm{B}}_{21}\cc_1\bar{\bm{B}}_{21}'+\tr\bar{\bm{B}}_{22}\cc_2\bar{\bm{B}}_{22}'\\
&=&\sum_{s=1}^2 
\tr(\hh^{(\covernum)}\g(\qq_{11},\qq_{22},\bm{0}))_{ss}\ddh_s\cc_s\
\ddh_s(\hh^{(\covernum)}\g(\qq_{11},\qq_{22},\bm{0}))_{ss}'\xxxi^{-1}_s
+ \tr\qqq_{21}\cc_1\qqq_{21}',
\end{eqnarray*}
\begin{eqnarray*}
\tr(\bar{\bm{B}}\cc\bar{\bm{B}}')^2
&=&\tr(\cc\bar{\bm{B}}'\bar{\bm{B}})^2\\
&=&\tr\left(
\begin{array}{cc}
\cc_1(\bar{\bm{B}}_{11}'\bar{\bm{B}}_{11}+\bar{\bm{B}}_{21}'\bar{\bm{B}}_{21})&
\cc_1\bar{\bm{B}}_{21}'\bar{\bm{B}}_{22}\\
\cc_2\bar{\bm{B}}_{22}'\bar{\bm{B}}_{21}& \cc_2\bar{\bm{B}}_{22}'\bar{\bm{B}}_{22}
\end{array}
\right)^2\\
&=&\tr(\cc_1(\bar{\bm{B}}_{11}'\bar{\bm{B}}_{11}+\bar{\bm{B}}_{21}'\bar{\bm{B}}_{21}))^2
+2\tr\cc_1\bar{\bm{B}}_{21}'\bar{\bm{B}}_{22}\cc_2\bar{\bm{B}}_{22}'\bar{\bm{B}}_{21}
+\tr(\cc_2\bar{\bm{B}}_{22}'\bar{\bm{B}}_{22})^2\\
&=&\tr\big(\cc_1\ddh_1(\hh^{(\covernum)}\g(\qq_{11},\qq_{22},\bm{0}))_{11}'
\xxxi_1^{-1}(\hh^{(\covernum)}\g(\qq_{11},\qq_{22},\bm{0}))_{11}\ddh_1
+\cc_1\qqq_{21}'\qqq_{21}\bigr)^2\\
& &+2\tr\big(\cc_1\qqq_{21}'\xxxi_2^{-1/2}(\hh^{(\covernum)}\g(\qq_{11},\qq_{22},\bm{0}))_{22}\ddh_2
\cc_2\\
& &\qquad\times\ddh_2(\hh^{(\covernum)}\g(\qq_{11},\qq_{22},\bm{0}))_{22}'\xxxi_2^{-1/2}\qqq_{21}
\big)\\
& &+\tr(\cc_2\ddh_2
(\hh^{(\covernum)}\g(\qq_{11},\qq_{22},\bm{0}))_{22}'\xxxi_2^{-1}(\hh^{(\covernum)}\g(\qq_{11},\qq_{22},\bm{0}))_{22}\ddh_2)^2.\\
\end{eqnarray*}
Consequently we have the following results; all the asymptotic expectations below are taken with respect to the distributions in (\ref{asympt_dist_lemma1}) and the spectral decompositions (\ref{decomp_W_ii}).
\begin{eqnarray}
\label{temp_form_E_tr}
\lefteqn{\lim_{\bdiva \rightarrow 0}E[\tr\g\Lllh\cc\Lllh\ggt\ggamma\Llambdain\ggammat]}\nonumber\\
&&=E[\tr\g_{11} \ddh_1 \cc_1\ddh_1 \ggt_{11} \xxxi_1^{-1}]
+E[\tr\g_{22} \ddh_2 \cc_2\ddh_2 \ggt_{22} \xxxi_2^{-1}]\\
& &\quad 
+E[\tr\bm{Z}_{21}\cc_1\bm{Z}_{21}']\nonumber
\label{temp_form_E_tr^2}
\end{eqnarray}
\begin{eqnarray}
\lefteqn{\lim_{\bdiva \rightarrow 0}E[\tr(\g\Lllh\cc\Lllh\ggt\ggamma\Llambdain\ggammat)^2]}\nonumber\\
&&=E[\tr(\cc_1\ddh_1 \ggt_{11} \xxxi_1^{-1}\g_{11} \ddh_1
+\cc_1\bm{Z}_{21}'\bm{Z}_{21})^2]\nonumber\\
& & \qquad
+2E[\tr\cc_1\bm{Z}_{21}'\xxxi_2^{-1/2}\g_{22} \ddh_2 \cc_2
\ddh_2 \ggt_{22} \xxxi_2^{-1/2}\bm{Z}_{21}]\nonumber\\
& &\qquad +E[\tr(\cc_2\ddh_2 \ggt_{22} \xxxi_2^{-1}\g_{22} \ddh_2 )^2]\nonumber\\
&&=E[\tr(\g_{11} \ddh_1 \cc_1\ddh_1 \ggt_{11} \xxxi_1^{-1})^2]\nonumber\\
&&\qquad +2\tr E[\cc_1\ddh_1 \ggt_{11} \xxxi_1^{-1}\g_{11} \ddh_1
\cc_1]E[\bm{Z}_{21}'\bm{Z}_{21}] \nonumber\\
& &\qquad +E[\tr\cc_1\bm{Z}_{21}'\bm{Z}_{21}\cc_1\bm{Z}_{21}'\bm{Z}_{21}]\\
& &\qquad +2\tr E[\xxxi_2^{-1/2}\g_{22} \ddh_2 \cc_2
\ddh_2 \ggt_{22} \xxxi_2^{-1/2}]E[\bm{Z}_{21}\cc_1\bm{Z}_{21}']\nonumber\\
& &\qquad +E[\tr(\g_{22} \ddh_2 \cc_2\ddh_2 \ggt_{22} \xxxi_2^{-1})^2].\nonumber
\end{eqnarray}

We further calculate the expectations related to $\bm{Z}_{21}$.
It is obvious that
\begin{equation}
\label{exp_rel_z}
E[\bm{Z}_{21}\cc_1\bm{Z}_{21}']=(\tr\cc_1)\ii_{p-m},\qquad E[\bm{Z}_{21}'\bm{Z}_{21}]=(p-m)\bm{I}_m,
\end{equation}
since $(\bm{Z}_{21})_{ij},\ 1\leq i \leq p-m,\ 1\leq j \leq m$, are all independently distributed as the standard normal distributions.
Letting $\bm{T}=(t_{ij})=\bm{Z}_{21}'\bm{Z}_{21}$ we have
\begin{eqnarray}
\label{exp_czzczz}
E[\tr\cc_1\bm{Z}_{21}'\bm{Z}_{21}\cc_1\bm{Z}_{21}'\bm{Z}_{21}]
&=&E[\tr\cc_1\bm{T}\cc_1\bm{T}]
=\sum_{i=1}^m c_iE[(\bm{T}\cc_1\bm{T})_{ii}]
\nonumber \\
&=&\sum_{i=1}^m c_iE[\sum_{s=1}^m t_{is}^2 c_s]
=\sum_{i=1}^m \sum_{s=1}^m c_ic_sE[t_{is}^2],
\end{eqnarray}
while
\begin{eqnarray}
\label{com_t_is}
E[t_{is}^2]
&=&E[(\sum_{j=1}^{p-m}(\bm{Z}_{21})_{ji}(\bm{Z}_{21})_{js})^2]\nonumber\\
&=&E[\sum_{j=1}^{p-m}(\bm{Z}_{21})^2_{ji}(\bm{Z}_{21})^2_{js}+2\sum_{j_1<j_2}(\bm{Z}_{21})_{j_1i}(\bm{Z}_{21})_{j_1s}(\bm{Z}_{21})_{j_2i}(\bm{Z}_{21})_{j_2s}].
\end{eqnarray}
We also have
\begin{equation}
\label{exp_2tuple_z}
E[(\bm{Z}_{21})_{ji}^2(\bm{Z}_{21})_{js}^2]=
\left\{
\begin{array}{cl}
3, & \mbox{ if }i=s,\\
1, & \mbox{ if }i\ne s,
\end{array}
\right.
\end{equation}
\begin{equation}
\label{exp_4tuple_z}
E[(\bm{Z}_{21})_{j_1i}(\bm{Z}_{21})_{j_1s}(\bm{Z}_{21})_{j_2i}(\bm{Z}_{21})_{j_2s}]=
\left\{
\begin{array}{cl}
1, & \mbox{ if }i=s,\\
0, & \mbox{ if }i\ne s.
\end{array}
\right.
\end{equation}
Substituting (\ref{exp_2tuple_z}),(\ref{exp_4tuple_z}) into (\ref{com_t_is}), we have
\begin{equation}
\label{cal_t_is^2}
E[t_{is}^2]=
\left\{
\begin{array}{cl}
(p-m)(p-m+2),& \mbox{ if }i=s,\\
 p-m,& \mbox{ if }i\ne s.
\end{array}
\right.
\end{equation}
Consequently from (\ref{exp_czzczz}) and (\ref{cal_t_is^2}),
\begin{equation}
\label{result_czzcc}
E[\tr\cc_1\bm{Z}_{21}'\bm{Z}_{21}\cc_1\bm{Z}_{21}'\bm{Z}_{21}]
=(p-m)(p-m+2)\sum_{i=1}^m c_i^2+2(p-m)\sum_{1\leq i<s \leq m}c_ic_s.
\end{equation}
Substituting (\ref{exp_rel_z}) and (\ref{result_czzcc}) into (\ref{temp_form_E_tr}), (\ref{temp_form_E_tr^2}), we have 
the result.
\subsection{Analytic Evaluation of Asymptotic Risk} 
\label{susec:Anal_Eval_Asymp_Risk}
We illustrate an analytic calculation of $E[d_i], E[d_i^{\:2}]$, 
$i=1,\ldots,p$, by the case $p=4,\ m=1$ and $n(\geq 4)$ even.
Suppose $\sss \sim \ww_3(n,\ii_p)$.
Note the density function of
$\lll=(l_1,l_2,l_3)$ is given by (see e.g. Theorem 3.2.18 of Muirhead
(1982))
$$
K_3(n)\prod_{i=1}^3 l_i^{\:\tn } (l_1-l_2)(l_1-l_3)(l_2-l_3)
\exp\Bigl(-\frac{1}{2}\sum_{i=1}^3l_i\Bigr),
$$
where $\tn = u(n)=(n-4)/2,$ which is an integer, and 
\begin{eqnarray*}
K_3(n)%
      &=&\pi^{3/2}/\Bigl(2^{3n/2}\Gamma(n/2)\Gamma((n-1)/2)\Gamma(n/2-1)\Gamma(3/2)\Gamma(1)\Gamma(1/2)\Bigr).
\end{eqnarray*}
Let 
$$
\Delta_1=l_1-l_2,\ \Delta_2=l_2-l_3,\ \Delta_3=l_3.
$$
The density function $f_3(\bm{\Delta})$ of $\bm{\Delta}=(\Delta_1,\Delta_2,\Delta_3)$ is given  by
\begin{eqnarray*}
f_3(\bm{\Delta})
&=&K_3(n)\Delta_3^{\tn }(\Delta_2+\Delta_3)^{\tn }(\Delta_1+\Delta_2+\Delta_3)^{\tn }\\
& &\quad \times\Delta_1\Delta_2(\Delta_1+\Delta_2)\exp\Bigl(-\frac{1}{2}(\Delta_1+2\Delta_2+3\Delta_3)\Bigr)\\
&=&K_3(n)\Bigl(\sum_{i=0}^{\tn } {\tn \choose i}%
\Delta_2^i\,\Delta_3^{\tn -i}\Bigr)
\Bigl(\sum_{s=0}^{\tn }\sum_{t=0}^{\tn -s}
{\tn \choose s} {\tn - s \choose t}
\: \Delta_1^s\,\Delta_2^t\,\Delta_3^{\tn -s-t}\Bigr)\\
& &\quad \times \Bigl(\sum_{j=0}^1\Delta_1^j\,\Delta_2^{1-j}\Bigr)\Delta_1\Delta_2\Delta_3^{\tn }\exp\Bigl(-\frac{1}{2}(\Delta_1+2\Delta_2+3\Delta_3)\Bigr)\\
&=&K_3(n) \sum_{i=0}^{\tn }\sum_{j=0}^1 \sum_{s=0}^{\tn }
\sum_{t=0}^{\tn -s}
{\tn \choose i} {\tn \choose s}{\tn -s \choose t}\;
\Delta_1^{j+s+1}\Delta_2^{i-j+t+2}\Delta_3^{3\tn -i-s-t}\\
& &\quad \times\exp\Bigl(-\frac{1}{2}(\Delta_1+2\Delta_2+3\Delta_3)\Bigr)
\end{eqnarray*}
We define a function $F_3(x_1,x_2,x_3;n)$ of nonnegative integers $x_i,\ i=1,2,3$, as
$$
F_3(x_1,x_2,x_3;n)=E[\Delta_1^{x_1}\Delta_2^{x_2}\Delta_3^{x_3}].
$$
Then
\begin{eqnarray*}
F_3(x_1,x_2,x_3;n)
&=&K_3(n) \sum_{i=0}^{\tn }\sum_{j=0}^1 \sum_{s=0}^{\tn }
\sum_{t=0}^{\tn -s}  
{\tn \choose i}{\tn \choose s} {\tn - s \choose t} \nonumber\\
& &\hspace{2cm}\times \int_0^\infty \Delta_1^{j+s+x_1+1}\exp\Bigl(-\frac{1}{2}\Delta_1\Bigr)d\Delta_1\nonumber\\
& &\hspace{2cm}\times \int_0^\infty \Delta_2^{i-j+t+x_2+2}\exp\Bigl(-\Delta_2\Bigr)d\Delta_2 \nonumber\\
& &\hspace{2cm}\times \int_0^\infty \Delta_3^{3\tn -i-s-t+x_3}\exp\Bigl(-\frac{3}{2}\Delta_3\Bigr)d\Delta_3\nonumber\\
&=&K_3(n) \sum_{i=0}^{\tn }\sum_{j=0}^1 \sum_{s=0}^{\tn }
\sum_{t=0}^{\tn -s}
{\tn \choose i}{\tn \choose s} {\tn - s \choose t} \nonumber\\
& &\hspace{2cm}\times 2^{3\tn -i+j-t+x_1+x_3+3}\:3^{-3\tn +i+s+t-x_3-1}\\
& &\hspace{2cm}\times (j+s+x_1+1)! (i-j+t+2+x_2)! \nonumber\\
& &\hspace{2cm}\times (3\tn -i-s-t+x_3)!\nonumber
\end{eqnarray*}

Note that for the case $p=4,\ m=1$, the distributions of $d_i$,
$i=1,\ldots,4$, in Theorem \ref{asym_risk_I&I} is given as follows;
$d_1=\ww_{11}\sim \chi^2_n$ and $d_2>d_3>d_4$ are the ordered
eigenvalues of $\ww_{22}\sim \ww_3(n-1,\ii_3).$ Using
$\Delta_1=d_2-d_3,\ \Delta_2=d_3-d_4,\ \Delta_3=d_4$ and
$F_3(x_1,x_2,x_3;n)$ as above, we can calculate
$\bm{b}=(b_1,\ldots,b_4)$ and $\bm{A}=(a_{ij})_{1\leq i,j \leq 4}$ in
Theorem \ref{asym_risk_I&I} as follows;
\begin{eqnarray*}
b_1&=&E[d_1]+(p-m)=n+3,\\
b_2&=&E[d_2]=E[\Delta_1+\Delta_2+\Delta_3]=F_3(1,0,0;n-1)+F_3(0,1,0;n-1)\\
& &\hspace{5.2cm}+F_3(0,0,1;n-1),\\
b_3&=&E[d_3]=E[\Delta_2+\Delta_3]=F_3(0,1,0;n-1)+F_3(0,0,1;n-1),\\
b_4&=&E[d_4]=E[\Delta_3]=F_3(0,0,1;n-1),\\
a_{11}&=&E[d_1^2+2(p-m)d_1]+(p-m)(p-m+2)\\
      &=&n^2+2n+6n+15=n^2+8n+15,\\
a_{22}&=&E[d_2^2]=E[\sum_{i=1}^3\Delta_i^2+2\sum_{1\leq i<j \leq 3}\Delta_i\Delta_j]\\
      &=&F_3(2,0,0;n-1)+F_3(0,2,0;n-1)+F_3(0,0,2;n-1)\\
      & &\;+2F_3(1,1,0;n-1)+2F_3(1,0,1;n-1)+2F_3(0,1,1;n-1),\\
a_{33}&=&E[d_3^2]=E[\Delta_2^2+\Delta_3^2+2\Delta_2\Delta_3]\\
      &=&F_3(0,2,0;n-1)+F_3(0,0,2;n-1)+2F_3(0,1,1;n-1),\\
a_{44}&=&E[d_4^2]=E[\Delta_3^2]=F(0,0,2;n-1),\\
a_{12}&=&a_{21}=E[d_2]=F_3(1,0,0;n-1)+F_3(0,1,0;n-1)+F_3(0,0,1;n-1),\\
a_{13}&=&a_{31}=E[d_3]=F_3(0,1,0;n-1)+F_3(0,0,1;n-1),\\
a_{14}&=&a_{41}=E[d_4]=F_3(0,0,1;n-1),\\
a_{23}&=&a_{32}=a_{24}=a_{42}=a_{34}=a_{43}=0.
\end{eqnarray*}

\end{document}